%%% Updated as of May 19th

%%%%%%%%%%%%%Abbreviations%%%%%%%%%%%%%%%%%%%%%

\font\gothic=eufm10.
\font\piccolo=cmr9.
\font\sets=msbm10.
.
\font\stampatello=cmcsc10.
.

\def\C{\hbox{\sets C}}

\def\N{\hbox{\sets N}}

\def\Z{\hbox{\sets Z}}
\def\defineq{\buildrel{def}\over{=}}

\def\doublesum{\mathop{\sum\sum}}

\def\SingSer{\hbox{\gothic S}}

\def\1{{\bf 1}}
\def\sgn{{\rm sgn}}
\def\square{\hbox{\vrule\vbox{\hrule\phantom{s}\hrule}\vrule}}

\centerline{\bf FINITE RAMANUJAN EXPANSIONS AND SHIFTED CONVOLUTION}
\centerline{\bf SUMS OF ARITHMETICAL FUNCTIONS, II}

\bigskip

\par
\centerline{\stampatello Giovanni Coppola and M. Ram Murty}

\bigskip
\bigskip
\bigskip

\footnote{}{
%\subjclass
$2010$ {\it Mathematics Subject Classification}: $11$A$25$,$11$K$65$,$11$N$37$

%\keywords
{\it Key words and phrases}: finite Ramanujan expansions, shifted convolution sum 

%\thanks
Research of the second author was partially supported by an NSERC Discovery grant. 
}

\par
\noindent
{\stampatello Abstract.}
{\piccolo 
We continue our study of convolution sums of two arithmetical functions $f$ and $g$, 
of the form $\sum_{n \le N} f(n) g(n+h)$, in the context of heuristic asymptotic formul\ae. 
Here, the integer $h\ge 0$ is called, as usual, the {\it shift} of the convolution sum.  
We deepen the study of finite Ramanujan expansions of general $f,g$ for the 
purpose of studying their convolution sum. Also, we introduce another kind of Ramanujan expansion 
for the convolution sum of $f$ and $g$, namely in terms of its shift $h$ and we compare this 
\lq \lq shifted Ramanujan expansion\rq \rq, with our previous finite expansions in terms of the $f$ and $g$ arguments. 
Last but not least, we give examples of such shift expansions, in classical literature, for the heuristic formul\ae. 
}

\bigskip
\bigskip

\par
\centerline{\stampatello 1. Introduction and statement of main results}%\label{Intro}1 
\bigskip
\par
\noindent
We start, as in our previous paper, from the definition of the {\it Ramanujan sum} (see [R] and compare [M] for the properties) : 
$$%\label{RamaSum}
c_q(n)\defineq \sum_{{a=1}\atop {(a,q)=1}}^{q}\cos\left({{2\pi an}\over q}\right) 
= \sum_{{a=1}\atop {(a,q)=1}}^{q}e^{2\pi ian/q} 
= \sum_{{d|q}\atop {d|n}}d\mu\left({q\over d}\right) 
\leqno{(1)}
$$
\par
\noindent
(compare first three eq.s in [M]), where we abbreviate with $(a,q)$ the greatest common divisor of any integers $a$ and $q$, as usual, with $\mu$ the 
M\"obius function (on primes 
$\mu(p)\defineq -1,$ and
$$  \mu(1)\defineq  1  , \qquad \mu(p_1 \cdots p_r)\defineq (-1)^r$$
 for $r\ge 2$ distinct primes $p_j$ and $\mu(n)\defineq 0$ on all other integers $n>1$). 
\par
Given $f,g:\N \rightarrow \C$ any arithmetic functions, we may consider, say, the {\it shifted convolution sum} of $f$ and $g$, which we abbreviate as the {\it correlation} of $f$ and $g$ (in the sequel), that we studied in our previous papers (in this series) 
$$%\label{CorDef}*
C_{f,g}(N,h)\defineq \sum_{n\le N}f(n)g(n+h),
\leqno{(2)}
$$
\par
\noindent
where the integer \enspace $h\ge 0$ \enspace is called the {\it shift}. Under suitable
conditions, we proved in [CMS2] that 
$$%\label{CorAsy}
C_{f,g}(N,h)=\SingSer_{f,g}(h)N+O(N^{1-\delta}(\log N)^{4-2\delta}) 
\leqno{(3)}
$$
\par
\noindent
(compare Theorem 2 [CMS2] for the precise statement), for a $\delta>0$, defining the {\it singular series} of $f$ and $g$ as 
$$%\label{defSingSer}*
\SingSer_{f,g}(h)\defineq \sum_{q=1}^{\infty}\widehat{f}(q)\widehat{g}(q)c_q(h), 
\leqno{(4)}
$$
\par
\noindent
where $\widehat{f}(q)$ and $\widehat{g}(q)$
are related to the {\it Ramanujan coefficients} of $f$ and $g$
respectively (see section $3$ for the definition and subsequent sections for the properties and examples).  
\par				% PAGE 2 
However, we briefly give the {\it Ramanujan expansion} of any $f$, of coefficients $\widehat{f}(q)$, 
$$%\label{RamaExp}
f(n)\defineq \sum_{q=1}^{\infty}\widehat{f}(q)c_q(n), 
\leqno{(5)}
$$
\par
\noindent
only assuming the pointwise convergence (compare Definition 2 in [CMS]). 
\par
Here, we recall the vital remark we made in [CMS2], in order to get, for fairly general $f$ and $g$, {\it finite} Ramanujan expansions (namely,  series like in $(5)$ become sums), defining for $f:\N \rightarrow \C$ the {\it Eratosthenes transform} (Aurel Wintner [W] coined this terminology), namely $f'\defineq f\ast \mu$, so M\"obius inversion [CojM] gives : $f(n)=\sum_{d|n}f'(d)$ (likewise for $g$) and we have 
$$%\label{VitaRema}
C_{f,g}(N,h)=\sum_{d}f'(d)\sum_{q}g'(q)\sum_{{n\le N}\atop {{n\equiv 0\, \bmod d}\atop {n\equiv -h\, \bmod q}}}1 
=\sum_{d\le N}f'(d)\sum_{q\le N+h}g'(q)\sum_{{n\le N}\atop {{n\equiv 0\, \bmod d}\atop {n\equiv -h\, \bmod q}}}1. 
\leqno{(6)}
$$
\par
\noindent
(It is suggestive to think of $f'$ as the ``arithmetical derivative'' of $f$.)
The above expression amounts to writing our (arbitrary) $f,g$ as {\it truncated divisor sums} (see next section's $(11)$, in turn, giving their finite Ramanujan expansions, see $(12)$ in $\S 3$). 
\par
We introduce, now, another possible approach, in the study of $f,g$ correlations. 
\par
In fact, apart from these {\sl finite} expansions (even if depending on $N,h$) that we introduced in [CMS2] (see $\S 3$), which are relative to the single (and arbitrary) functions $f$ and $g$, we may consider the {\stampatello shift-Ramanujan expansion} (or \lq \lq Ramanujan expansion with respect
to the shift\rq \rq), abbreviated { s.R.e.}, of $f,g$ correlation, namely 
$$%\label{sRe}
C_{f,g}(N,h)=\sum_{\ell=1}^{\infty}\widehat{C_{f,g}}(N,\ell)c_{\ell}(h)
\leqno{(7)}
$$
\par
\noindent
where, now, the main issue is the possibility to give such an expansion, with some \lq \lq {\stampatello shift-Ramanujan coefficients}\rq \rq, $\widehat{C_{f,g}}(N,\ell)$, and whether we have in $(7)$ an absolutely or uniformly convergent series. Many classical results in the literature, like our results above for $C_{f,g}(N,h)$, are all pointing towards the heuristic formula for these coefficients 
$$%\label{CfgHat}
\widehat{C_{f,g}}(N,\ell)
\sim
\widehat{f}(\ell)\widehat{g}(\ell)N
\leqno{(8)}
$$
\par
\noindent
where the \lq \lq $\sim$\rq \rq \enspace sign is used like for Fourier coefficients formul\ae, i.e., after suitable analytic assumptions (and, also, with a well-specified analytic meaning). 
\bigskip
\par
\noindent
The analytic assumptions ensuring \lq \lq good\rq \rq \thinspace convergence may be very complicated, for these shift-Ramanujan expansions. However, the above for s.R.e. coefficients look like well-known heuristic formul\ae, starting with the Hardy-Littlewood conjecture on prime tuples [HLi] (see $\S5$ end). 
\par
In the following, \lq \lq $f$ is essentially bounded\rq \rq, i.e. $f(n)\ll_{\varepsilon} n^{\varepsilon}$, is tantamount to \lq \lq $f$ satisfies the Ramanujan Conjecture\rq \rq. Hereafter, \lq \lq $\forall \varepsilon>0$\rq \rq, as usual, is implicit in bounds; in fact, $\varepsilon>0$ is arbitrarily small and may change even in the same formula. 
\par
We give, inspired by these heuristics, the following general bounds, for all real $\delta>0$ : we say, by definition, that a \lq \lq s.R.e. is in $\delta-$class\rq \rq, whenever (for essentially bounded $f,g$) 
$$%\label{deltaCl}
\widehat{C_{f,g}}(N,\ell)\ll {{N^{1+\varepsilon}}\over {\ell^{1+\delta}}}, 
\leqno{(9)}
$$
\par
\noindent
with the implied constant depending eventually on both $\delta,\varepsilon$; the
noteworthy case $\delta=1$ will be referred to as \lq \lq s.R.e. is in the 
first class\rq \rq. For example, equation $(8)$ above implies (by the bounds on $f,g$ Ramanujan coefficients of $(14)$ in $\S 3$) that our s.R.e. is in the first class (for the essentially bounded $f,g$ and shift $h\ll N$, assuming also remainders in $(8)$ are small enough). 
\par				% PAGE 3 
We define, here, a {\it pure} Ramanujan expansion by \enspace $F(v)=\sum_{q=1}^{\infty}\widehat{F}(q)c_q(v)$, \enspace (pointwise converging in $v$ and) in which the $v-$dependence is only in $c_q(v)$. For example, take [M], p.24, in which both $\widehat{F}(q)=1/q$ and $\widehat{F}(q)=0$ (on all $q$) represent the  constant zero function. On the other hand, Hildebrand's finite Ramanujan expansions $(1.4)$, p.167 [SchSpi] aren't pure in our sense. 
\smallskip
\par
Our main result is the following. 
\enspace Recall {\it Euler's function} is \thinspace $\varphi(\ell)\defineq |\{ n\le \ell : (n,\ell)=1 \}|$. 
\smallskip
\par
\noindent {\bf Theorem 1.} {\it Let } $N,h\in \N$, {\it with the shift } $h\ll N$, {\it as } $N\to \infty$, {\it and assume that } $f,g:\N \rightarrow \C$ {\it are essentially bounded, with, say, } $\max\{q\ge 1:g'(q)\neq 0\}$
{\it as a bound and } $f,g$ {\it not depending on } $h$. {\it Consider the shift-Ramanujan expansion } $(7)$, {\it abbreviated s.R.e., assuming, as we may\footnote{$^1$}{\rm by Hildebrand's Theorem, see $\S1$ end in [M] and reference [10] therein, also, compare beginning of $\S 5$.}, that it converges pointwise, for all the fixed} \thinspace $h\in \N$. ({\it We don't assume uniqueness of this expansion: we may even have undetermined coefficients}.). {\it Then, the following are equivalent. } 
\par
\noindent
  \item{$(i)$} {\it the s.R.e. is uniformly convergent} ({\it i.e., $(24)$ in} $\S 5$) {\it and pure};
  \item{$(ii)$} {\it the s.R.e. coefficients are given by Carmichael's formula}
$$
\widehat{C_{f,g}}(N,\ell)={1\over {\varphi(\ell)}}\lim_{x\to \infty}{1\over x}\sum_{h\le x}C_{f,g}(N,h)c_{\ell}(h); 
$$
  \item{$(iii)$} {\it the s.R.e. coefficients are given by the explicit formula}
$$
\widehat{C_{f,g}}(N,\ell)={{\widehat{g}(\ell)}\over {\varphi(\ell)}}\sum_{n\le N}f(n)c_{\ell}(n); 
$$
  \item{$(iv)$} {\it the s.R.e. is finite and pure}. 
{\it In the same hypotheses, whenever any one of these conditions holds, the s.R.e. is in the first class.} 
\smallskip
\par
\noindent {\bf Remark.} The \thinspace $\max\{q\ge 1:g'(q)\neq 0\}$ \thinspace exists, since $(11)$ in $\S 2$ gives $q>N+h \Rightarrow g'(q)=0$. 

\smallskip

\par
We will call a s.R.e. satisfying one (hence, all) of previous equivalent conditions a {\it regular} shift-Ramanujan expansion. 
\par
An example of regular s.R.e. is the main term in $(3)$, since 
$\widehat{f}(q)=0$, $\forall q>N$, by $(12)$ (compare next section's final remarks on the singular series as a sum). 

However, there are examples of irregular s.R.e., one of which (an arithmetic function $f_H$ depending on the parameter $H$) is given in $\S 9$. 

\medskip

\par
We prove in $\S 5$ the following immediate, very important consequence of our Theorem 1. 
\smallskip
\par
\noindent {\bf Corollary 1.} {\it With the hypotheses of Theorem 1, when the s.R.e. is regular and the $f$ Eratosthenes transform} $f'$ {\it is supported up to $D$} ({\it i.e.,} $f'(d)=0, \forall d>D$), {\it with $D\ll N$ as $N\to \infty$, we get } 
$$
C_{f,g}(N,h)=\SingSer_{f,g}(h)N+O_{\varepsilon}\left(N^{\varepsilon}D\right). 
$$
\par
\noindent
{\it In particular, when } ${{\log D}\over {\log N}}<1$, {\it say} \thinspace ${{\log D}\over {\log N}}<1-\delta$ \thinspace {\it for a certain small $\delta>0$, we have}
$$
C_{f,g}(N,h)=\SingSer_{f,g}(h)N+O\left(N^{1-\delta + \epsilon}\right). 
$$

\medskip

\par
\noindent
We explicitly remark that here we have given our main results, but along the paper (especially in the last section) we will give other \lq \lq minor\rq \rq, so to  speak, results about these arguments. 

\vfill
\eject

\par				% PAGE 4 
The paper is organized as follows: 
\medskip
\par
\noindent
  \item{$\diamondsuit$} in the next section we highlight links between correlations of two (arbitrary) arithmetic functions $f$ \& $g$ and $f$ \& $g$ truncated divisor sums counterparts; 
  \item{$\diamondsuit$} in $\S 3$ we show the duality between those truncated divisor sums and their finite Ramanujan expansions with their properties; 
  \item{$\diamondsuit$} in $\S 4$ we give some examples of finite Ramanujan expansions; 
  \item{$\diamondsuit$} in $\S 5$ the new \lq \lq shift Ramanujan expansions\rq \rq, for correlations of $f$ and $g$, show their links with previous finite expansions for single $f$ and $g$, especially in the two proofs, of Theorem 1 \& Corollary 1; 
  \item{$\diamondsuit$} we make a short \lq \lq detour\rq \rq, regarding sieve functions in our context, in $\S 6$; 
  \item{$\diamondsuit$} we show how finite Ramanujan expansions change, in case our arithmetic functions, say, \lq \lq have no small prime divisors\rq \rq, in $\S 7$; 
  \item{$\diamondsuit$} in $\S 8$ we make some further remarks; 
  \item{$\diamondsuit$} finally in Appendix $\S 9$ we give useful lemmas
of independent interest. 

\bigskip
\bigskip

\par
\centerline{\stampatello 2. Shifted convolution sums and truncated divisor sums}%\label{SCStds}2
\bigskip
\par
\noindent
The heuristics for our correlations (compare classic [HLi]: eq. (5.26) and Conjecture B, with papers [CMS,CMS2,GMP,MS]) are of the kind 
$$%\label{heurCorr}
C_{f,g}(N,h)=\SingSer_{f,g}(h)N+\hbox{\rm good\enspace remainder}
\leqno{(10)} 
$$
\par
\noindent
(say, $O(N^{1-\delta})$, for a $\delta>0$), the singular series $\SingSer_{f,g}(h)$ for $f$ and $g$, of shift $h\ge 0$, being defined as above. 
A justification for this heuristic comes from the following considerations.
\par
For any $f,g:\N \rightarrow \C$ we defined Eratosthenes transforms $f'$ and $g'$, so $f(n)=\sum_{d|n}f'(d)$ and $g(m)=\sum_{q|m}g'(q)$; from our previous paper [CMS2] we know that, for our purposes (namely, confining to $C_{f,g}$ study), $f,g$ may be truncated over the divisors, as in $(6)$, getting 
$$%\label{fgTrunc} (11)
f(n)=\sum_{d|n,d\le N}f'(d)
\qquad \hbox{\rm and} \qquad
g(m)=\sum_{q|m,q\le N+h}g'(q).
\leqno{(11)}
$$
\par
\noindent
Thus in studying $C_{f,g}$ we naturally find the {\it finite Ramanujan expansions} of $f$ and $g$ ( $(12)$ in the next section). We use these truncated divisor sum representations for $f$ and $g$ to deduce: 
$$
C_{f,g}(N,h)=\sum_{d\le N}f'(d)\sum_{q\le N+h}g'(q)\sum_{{n\le N}\atop {{n\equiv 0\, \bmod d}\atop {n\equiv -h\, \bmod q}}}1 
$$
$$
=\sum_{l|h}\sum_{d\le N}f'(d)\sum_{{q\le N+h}\atop {(q,d)=l}}g'(q)\sum_{{n\le N}\atop {{n\equiv 0\, \bmod d}\atop {n\equiv -h\, \bmod q}}}1 
=\sum_{{l|h}\atop {b:=-h/l}}\sum_{t\le {N\over l}}f'(lt)\sum_{{r\le {{N+h}\over l}}\atop {(r,t)=1}}g'(lr)\sum_{{m\le {N\over {lt}}}\atop {m\equiv \overline{t}b\, \bmod r}}1, 
$$
\par
\noindent
defining \thinspace $\overline{t}\bmod r$  with $\overline{t}t\equiv 1\bmod r$. Whence an approach based on writing 
$$
\sum_{{m\le {N\over {lt}}}\atop {m\equiv \overline{t}b\, \bmod r}}1\sim {N\over {ltr}} = {{Nl}\over {(lt)(lr)}} 
$$
\par
\noindent
(meaning, of course, that fractional parts are assumed negligible here) gives the heuristic (again, $\sim$ has to be given the right meaning) 
$$
C_{f,g}(N,h)\sim \sum_{l|h}\sum_{d\le N}f'(d)\sum_{{q\le N+h}\atop {(q,d)=l}}g'(q){{Nl}\over {dq}} 
=N\sum_{l|h}l\sum_{d\le N}{{f'(d)}\over d}\sum_{{q\le N+h}\atop {(q,d)=l}}{{g'(q)}\over q} 
$$
\par				% PAGE 5 
\noindent
which is, recalling that we are truncating $f'$ at $N$ and $g'$ at $N+h$, 
$$
C_{f,g}(N,h)\sim N\sum_{l|h}l\sum_{d}{{f'(d)}\over d}\sum_{{q}\atop {(q,d)=l}}{{g'(q)}\over q} 
$$
\par
\noindent
which is the heuristic $(10)$, since we prove in $\S 9$, Lemma A.6, that 
$$
\SingSer_{f,g}(h)=\sum_{l|h}l\sum_{d}{{f'(d)}\over d}\sum_{{q}\atop {(q,d)=l}}{{g'(q)}\over q}. 
$$

\bigskip

\par
In these formul\ae \enspace for the singular series, we applied the absolute convergence of each series involved: all of them, thanks to the fact that our Ramanujan expansions are {\it finite}, and so clearly converge absolutely, hence justifying, in our context, all the series exchanges (compare Lemma A.6 proof, in $\S 9$). The same singular series, because of this fact, is simply a \lq \lq {\it singular sum}\rq \rq. This feature (like, also, our considerations on finiteness of Ramanujan expansions involved) seems to have been overlooked
 in the literature. Actually, our {\it singular sum} may be seen as the $N-$th partial sum, of the original singular series. We leave, as an exercise for the interested reader, to prove that the tail (as usual, the difference between the whole series and the partial sum) converges (very rapidly indeed) to zero (with a dependence, of course on $N$). 

\bigskip
\bigskip

\par
\centerline{\stampatello 3. The finite Ramanujan expansions: properties and formul\ae}%\label{fRe}3
\bigskip
\par
\noindent
The truncated divisor sums for $f$ and $g$, in $(11)$, have (compare [CMS2] Introduction) finite Ramanujan expansions, that we will  sometimes abbreviate { f.R.e.}: 
$$%\label{(0)} (12)
f(n)=\sum_{r\le N}\widehat{f}(r)c_r(n)
\qquad \hbox{\rm and} \qquad 
g(m)=\sum_{s\le N+h}\widehat{g}(s)c_s(m), 
\leqno{(12)}
$$
\par
\noindent
with an {\it explicit formula for} their {\it Ramanujan coefficients} that we proved in [CMS2] Introduction (soon after Lemma 1): 
$$%\label{(1)} (13)
\widehat{f}(r)=\sum_{{m}\atop {m\equiv 0\bmod r}}{{f'(m)}\over m}
={1\over r}\sum_{n\le {N\over r}}{{f'(rn)}\over n},
\quad
\widehat{g}(s)=\sum_{{m}\atop {m\equiv 0\bmod s}}{{g'(m)}\over m}
={1\over s}\sum_{n\le {{N+h}\over s}}{{g'(sn)}\over n}. 
\leqno{(13)}
$$
\par
\noindent
Notice that this formula implies all truncated divisor sums have finite Ramanujan expansion. 
\par
In particular, for the essentially bounded $f,g$, we get the bounds 
$$%\label{fReBounds} (14)
\widehat{f}(r)\ll_{\varepsilon} {{N^{\varepsilon}}\over r},
\quad
\widehat{g}(s)\ll_{\varepsilon} {{(N+h)^{\varepsilon}}\over s}. 
\leqno{(14)}
$$
\par
\noindent
The other {\it explicit formula}, this time {\it for the Eratosthenes transform} in terms of Ramanujan coefficients (see the Introduction of [CMS2], soon before Theorem 1), is : 
$$%\label{(2)} (15)*
f'(d)=d\sum_{j=1}^{\infty}\mu(j)\widehat{f}(dj) 
=d\sum_{j\le {N\over d}}\mu(j)\widehat{f}(dj), 
\enspace 
g'(q)=q\sum_{j=1}^{\infty}\mu(j)\widehat{g}(qj) 
=q\sum_{j\le {{N+h}\over q}}\mu(j)\widehat{g}(qj). 
\leqno{(15)}
$$
\par
\noindent
We profit to prove it briefly: from the M\"{o}bius inversion formula \enspace $\sum_{j|K}\mu(j)=[1/K]$ \enspace (with $[\enspace ]$ the integer part:   see {CMS2}, Lemma 3), together with $(13)$, 
$$
d\sum_{j\le {N\over d}}\mu(j)\widehat{f}(dj)=\sum_{j\le {N\over d}}{{\mu(j)}\over j}\sum_{n\le {N\over {dj}}}{{f'(djn)}\over n}
=\sum_{K\le {N\over d}}{{f'(dK)}\over K}\sum_{j|K}\mu(j)
=f'(d). 
$$
\par				% PAGE 6 
\noindent
These formul\ae \enspace link Ramanujan coefficients (resp., $\widehat{f}$, $\widehat{g}$), with Eratosthenes transforms (resp., $f'$, $g'$). 
This is a kind of duality: truncated divisor sums (with $f',g'$) and finite Ramanujan expansions (with $\widehat{f},\widehat{g}$) clearly describe the same objects (our functions $f,g$). 

Furthermore, for the {\it high} Ramanujan coefficients, i.e., having index \enspace $Q/2<q\le Q$, when divisors are truncated at $Q$, i.e., Eratosthenes transform support $\subseteq [1,Q]$, we have : 
$$%\label{highRama} (16)*
u(n)=\sum_{d|n,d\le Q}u'(d)
\quad \Longrightarrow \quad
\widehat{u}(q)=u'(q)/q, 
\enspace 
\forall q\in(Q/2,Q], 
\leqno{(16)}
$$
\par
\noindent
entailing 
$$%\label{highRf,g} (17)
\widehat{f}(r)={{f'(r)}\over r}, 
\enspace 
\forall r\in \left({N\over 2},N\right] 
\quad \hbox{\rm and} \quad 
\widehat{g}(s)={{g'(s)}\over s}, 
\enspace 
\forall s\in \left({{N+h}\over 2},N+h\right]. 
\leqno{(17)}
$$

\bigskip
\bigskip

\par
\centerline{\stampatello 4. The finite Ramanujan expansions: examples}%\label{fReExs}4
\bigskip
\par
\noindent
These two formul\ae \enspace in $(17)$ immediately imply for the von Mangoldt function, say, 

$$
\Lambda_N(n)=\sum_{d|n,d\le N}(-\mu(d)\log d), 
$$
\par
\noindent
which has been truncated as above for the calculation of $C_{\Lambda,\Lambda}(N,h)$, 

$$%\label{highLambda} (18)
\widehat{\Lambda_N}(r)=-{{\mu(r)\log r}\over r}, 
\thinspace 
\forall r\in \left({N\over 2},N\right] , 
\widehat{\Lambda_{N+h}}(s)=-{{\mu(s)\log s}\over s}, 
\thinspace 
\forall s\in \left({{N+h}\over 2},N+h\right]. 
\leqno{(18)}
$$
\par
\noindent
More precisely,  

$$
C_{\Lambda,\Lambda}(N,h)=\sum_{n\le N}\Lambda(n)\Lambda(n+h)
=\sum_{d\le N}\mu(d)(\log d)\sum_{q\le N+h}\mu(q)(\log q)\sum_{{n\le N}\atop {{n\equiv 0\bmod d}\atop {n+h\equiv 0\bmod q}}}1
$$
$$
=\sum_{d\le N}\mu(d)(\log d)\sum_{q\le N}\mu(q)(\log q)\sum_{{n\le N}\atop {{n\equiv 0\bmod d}\atop {n+h\equiv 0\bmod q}}}1
 -\sum_{N<q\le N+h}\mu(q)(\log q)\sum_{{n\le N}\atop {n+h\equiv 0\bmod q}}\Lambda(n) 
$$
$$
=\sum_{d\le N}\mu(d)(\log d)\sum_{q\le N}\mu(q)(\log q)\sum_{{n\le N}\atop {{n\equiv 0\bmod d}\atop {n+h\equiv 0\bmod q}}}1
 +O\left(\log^2(N+h)\sum_{N<q\le N+h}\left({N\over q}+1\right)\right). 
$$
\par
\noindent
Remainder terms are clearly $\ll h\log^2(N+h)$, which, if $h=o(N)$ is small enough, say $h\ll N^{1-\delta}$ (for a fixed $\delta>0$), are negligible. 
\par
Hence we may stick to only one truncation, the one with $N$ (ignoring the shift). Of course, for $h$ small enough, this procedure works for all arithmetic functions $f,g$ that do not grow too fast, like the ones satisfying (like $\Lambda$) the Ramanujan Conjecture: $f(n),g(n)\ll n^{\varepsilon}$. 
\par
Since all of our examples will satisfy this growth condition (like all interesting arithmetic functions; otherwise we may re-normalize) we can always write, in good approximation: 
$$%\label{easyCfg} (19)*
C_{f,g}(N,h)=\sum_{d\le N}f'(d)\sum_{q\le N}g'(q)\sum_{{n\le N}\atop {{n\equiv 0\bmod d}\atop {n+h\equiv 0\bmod q}}}1, 
\leqno{(19)}
$$
\par				% PAGE 7 
\noindent
i.e., set a common truncation for both $f$ and $g$: 
$$ 
f(n)=\sum_{d|n, d\le N}f'(d)
\qquad \hbox{\rm and} \qquad 
g(m)=\sum_{q|m, q\le N}g'(q). 
$$
\par
\noindent
Returning to our von Mangoldt function $f=g=\Lambda$ (hence, $C_{f,g}(N,h)$ regards $h$-twins of primes), the idea of truncating its divisors has been pursued by many authors in the literature (mainly, in the area of sieve methods) and, very recently, has given spectacular results at the hands of Goldston, Pintz \& Y{\i}ld{\i}r{\i}m in the 2000s which had been applied by Green \& Tao (to prove that the sequence of primes contains arbitrary long arithmetic progressions), by Zhang, Maynard \& Tao, Polymath project to study bounded gaps between primes. 

We wish to emphasize that such an approach has not yet been followed, in order to give hints in the asymptotic formul\ae, for the correlation sum of twin primes (say, $C_{\Lambda,\Lambda}$ here) ! Thus, with $(18)$ above, we try to give a new flavor to the estimate of Ramanujan coefficients of $\Lambda$; these high coefficients are somehow unexpected as they do not agree, exactly, with the known classical ones.
\par
However, our formul\ae, specifically $(13)$, give for the low coefficients (from now on we will work with $N-$truncations) 
$$%\label{lowCoef} (20)
\widehat{\Lambda}(q)=-{1\over q}\sum_{n\le {N\over q}}{{\mu(qn)\log(qn)}\over n}
=-{{\mu(q)\log q}\over q}\sum_{{n\le {N\over q}}\atop {(n,q)=1}}{{\mu(n)}\over n}
 -{{\mu(q)}\over q}\sum_{{n\le {N\over q}}\atop {(n,q)=1}}{{\mu(n)\log n}\over n}
\leqno{(20)}
$$
\par
\noindent
which, heuristically speaking (for $q$ small with respect to $N$), are in good agreement with the classical known formul\ae, i.e. 
$$
\sim {{\mu(q)}\over {\varphi(q)}}, 
$$
\par
\noindent
from the very well-known formul\ae, see [CojM] or [MoV], for $q\le x$ (so, for $q\le \sqrt{N}$ in our case), $c>0$ fixed: 
$$%\label{muApprox} (21)
\sum_{{n\le x}\atop {(n,q)=1}}{{\mu(n)}\over n}\ll \exp\left(-c\sqrt{\log x}\right),
\quad 
\sum_{{n\le x}\atop {(n,q)=1}}{{\mu(n)\log n}\over n}=-{q\over {\varphi(q)}}+O\left(\exp\left(-c\sqrt{\log x}\right)\right). 
\leqno{(21)}
$$
\par
\noindent
The fact that we are  working with f.R.e. (equivalently, of truncating divisors) makes the coefficients behave in a different way, with respect to the classical Ramanujan (series) expansions; in particular, the \lq \lq low\rq \rq coefficients (i.e., with $q\le \sqrt{N}$ in above example) should agree (to some extent) with the classical ones, while the effect of truncating divisors is clear on the last ones (we call \lq \lq high\rq \rq, see the above), which may be totally different !

\bigskip

\par
We consider, now, three other examples, the first two of which are related to [CMS], respectively to Corollary 1 of [CMS] and to Corollary 2 of [CMS]. Our third of these (and last example) will be related to these two (but we didn't mention it, in our earlier papers). 

\medskip

\par
Our next example comes from the arithmetic function, say (see Corollary 1 [CMS]), 
$$
f_s(n)={{\sigma_s(n)}\over {n^s}}
={1\over {n^s}}\sum_{d|n}d^s
=\sum_{d|n}\left({n\over d}\right)^{-s}
=\sum_{d|n}d^{-s}
=\sigma_{-s}(n), 
$$
\par
\noindent
by passing from $d|n$ to its complementary divisor \thinspace $\left.{n\over d}\right|n$ \thinspace (we'll refer to this trick of Dirichlet as \lq \lq {\it flipping the divisors}\rq \rq). For this function we have, for all $s>0$, 
$$
\sigma_{-s}(n)={{\sigma_s(n)}\over {n^s}}
=\sum_{q=1}^{\infty}\widehat{\sigma_{-s}}(q)c_q(n) 
$$
\par				% PAGE 8 
\noindent
as a classical Ramanujan expansion (even {\it converging absolutely}, by Lemma 1 in $\S 9$). 
\par
Here for notational convenience we write \enspace $f_s(n):={{\sigma_s(n)}\over {n^s}}=\sigma_{-s}(n)$, introducing 
$$%\label{(joint)}22*
f_s(n)=\sum_{d|n}f'_s(d),
\quad {\rm with } \enspace f'_s\defineq f_s\ast \mu 
\quad \Longrightarrow \quad 
f_{s,D}(n)\defineq \sum_{{d|n}\atop {d\le D}}f'_s(d),
\leqno{(22)}
$$
\par
\noindent
(for fixed $s>0$ and $D\in \N$) its {\it truncated} counterpart, {\it over the divisors} up to $D$. This definition is pretty general; here, in the present example, 
$$
f'_s(d)=d^{-s}
\quad \Longrightarrow \quad 
\sigma_{-s,D}(n)\defineq \sum_{{d|n}\atop {d\le D}}d^{-s}
=\sum_{q\le D}\widehat{\sigma_{-s,D}}(q)c_q(n) 
$$
\par
\noindent
(for all $s>0$ and $D\in \N$, both fixed), with finite Ramanujan coefficients 
$$
\widehat{\sigma_{-s,D}}(q)\defineq \sum_{{m\le D}\atop {m\equiv 0\bmod q}}m^{-s-1}
={1\over {q^{s+1}}}\sum_{n\le {D\over q}}{1\over n^{s+1}}. 
$$
\par
\noindent
Notice that these are different from the classical Ramanujan coefficients, that we calculated thanks to the Delange 1976 Theorem (see Theorem 3 and following discussion, before Theorem 4, in [M]) 
$$
\widehat{\sigma_{-s}}(q)\defineq \sum_{m\equiv 0\bmod q}m^{-s-1}
={1\over {q^{s+1}}}\sum_{n=1}^{\infty}{1\over n^{s+1}}
={1\over {q^{s+1}}}\zeta(s+1) 
$$
\par
\noindent
(apart from similarity we'll check soon, for suitable indices),  since the  finite ones have
$$
{D\over 2}<q\le D 
\enspace \Rightarrow \enspace 
\widehat{\sigma_{-s,D}}(q)={1\over {q^{s+1}}}, 
$$
\par
\noindent
as we already knew from $(17)$, in the previous $\S 3$. (Also, these trivially vanish for $q>D$, while the classical ones don't.) 
\smallskip
\par 
However, if the indices are somehow \lq \lq small\rq \rq, the f.R.e. has coefficients that are asymptotic to classical ones, as we see now : when $D\to \infty$, 
$$
\widehat{\sigma_{-s,D}}(q)={1\over {q^{s+1}}}\sum_{n\le {D\over q}}{1\over n^{s+1}}
={\zeta(s+1)\over {q^{s+1}}}-{1\over {q^{s+1}}}\sum_{n>{D\over q}}{1\over n^{s+1}}
={\zeta(s+1)\over {q^{s+1}}}\left(1+O_s\left(\left({q\over D}\right)^s\right)\right), 
$$
\par
\noindent
namely for
$$
q=o(D)  
\enspace \Rightarrow \enspace 
\widehat{\sigma_{-s,D}}(q)=\widehat{\sigma_{-s}}(q)\left(1+o_s(1)\right)
\sim \widehat{\sigma_{-s}}(q).  
$$
\par
\noindent
Thus the finite Ramanujan coefficients, on \lq \lq small indices\rq \rq (say, $q=o(D)$ \thinspace here) are asymptotic to the Ramanujan coefficients (of the \lq \lq classical\rq \rq expansion). 

\medskip

\par
Furthermore, with \enspace $f_s(n)=\prod_{p|n}\left(1-p^{-s}\right)$, $s>0$, we get Corollary 2 [CMS] application and, from multiplicativity, 
$$
f_s(n)=\sum_{d|n}\mu(d)d^{-s}
\enspace \Rightarrow \enspace 
f'_s(d)=\mu(d)d^{-s}, 
$$
\par
\noindent
giving rise again to 
$$%\label{asympCoeFsD} (23)
q=o(D)  
\enspace \Rightarrow \enspace 
\widehat{f_{s,D}}(q)=\widehat{f_{s}}(q)\left(1+o_s(1)\right)
\sim \widehat{f_{s}}(q), 
\leqno{(23)}
$$
\par
\noindent
by the same calculations as above. 

\medskip

\par				% PAGE 9 
This suggests the third example of this kind, i.e., say  
$$
f_s(n)=\prod_{p|n}\left(1+p^{-s}\right)
=\sum_{d|n}\mu^2(d)d^{-s}
\enspace \Rightarrow \enspace 
f'_s(d)=\mu^2(d)d^{-s}
$$
\par
\noindent
with this kind of truncation 
$$
f_{s,D}(n)=\sum_{{d|n}\atop {d\le D}}\mu^2(d)d^{-s}
\enspace \Rightarrow \enspace 
f_{s,D}(n)=\sum_{q\le D}\widehat{f_{s,D}}(q)c_q(n) 
$$
\par
\noindent
with the same behavior as given in $(23)$. 

\medskip

\par
More generally Delange's Theorem gives (compare Theorem 3 on [M]) 
$$
\widehat{f}(q)=\sum_{m\equiv 0\bmod q}{{f'(m)}\over m}, 
$$
\par
\noindent
in the hypothesis 
$$
\sum_{m=1}^{\infty}{{2^{\omega(m)}|f'(m)|}\over m}<\infty, 
$$
\par
\noindent
which is certainly satisfied by the $D-$truncation of our $f$ (since we have a f.R.e. for it), say 
$$
f_D(n)=\sum_{{d|n}\atop {d\le D}}f'(d)
\enspace \Rightarrow \enspace 
\widehat{f_D}(q)=\sum_{{m\le D}\atop {m\equiv 0\bmod q}}{{f'(m)}\over m}, 
$$
\par
\noindent
which we proved directly (actually, the same method, but applying analytic approximations, too, of course, proves Delange's Theorem, compare [M]). 
We wish to have \enspace $\widehat{f_D}(q)\sim \widehat{f}(q)$ ! 
\par
For this, the only  problem is the effect of truncation on Ramanujan coefficients, that is, 
\bigskip
\par
\centerline{ differences ${\displaystyle \sum_{{m>D}\atop {m\equiv 0\bmod q}}{{f'(m)}\over m}={1\over q}\sum_{n>D/q}{{f'(qn)}\over n} }$  must be infinitesimal}
\bigskip
\par
\noindent
and this is possible, clearly, only when the variable \enspace $D/q\to \infty$, i.e. \thinspace $q=o(D)$ \thinspace here is a necessary condition (actually, {\it for previous three cases} a {\it sufficient} one, too). For the high coefficients we already observed a neat difference with the classical ones; this may be justified by the truncation (of divisors) itself, that has to change \lq \lq last\rq \rq, so to speak, coefficients, in order to cope with the infinite tail, that is missing (of course, in {\it finite} Ramanujan expansions). 

\bigskip
\bigskip

\par
\centerline{\stampatello 5. Ramanujan expansions with respect to the shift}%\label{RamaShift}5 
\bigskip
\par
\noindent
In the following, we will dwell mainly with the easiest possible hypothesis for the series in $(7)$, namely, the uniform convergence 
(i.e., not depending on the shift $h$) 
$$%\label{sReAbs} (24)
\sum_{\ell=1}^{\infty}\widehat{C_{f,g}}(N,\ell)c_{\ell}(h) 
\enspace
\hbox{\rm converges}, 
\enspace
\forall h\in \N.
\leqno{(24)}
$$
\par
\noindent
We explicitly point out that we are considering series back again, since our previous remark, truncating the divisor sums (hence giving f.R.e. above), gives no hint on whether the present expansion in $(24)$ is finite or not. 
\par
However, in the case of s.R.e. regularity then (see Theorem 1) the shift expansion itself is finite.  Our main problem is to try to understand {\sl when} we have such  regularity! 
\par				% PAGE 10 
Even in the case of irregularity, we always have the pointwise convergence of our s.R.e., as an easy consequence of Hildebrand's Theorem (see Theorem 1 footnote in $\S 1$) for any arithmetic function (here, the $f,g$ correlation as a function of the shift $h\in \N$); the problem, however, is that we don't have, a priori, the uniqueness for the s.R.e. and this may lead even to more different coefficients, for the same expansion. So, uniform convergence of our s.R.e. confirms to be the easiest analytic assumption, especially in the light of Theorem 1. 

\medskip

\par
We prove first Theorem 1 and, then, the much easier Corollary 1. 

\medskip

\par
\noindent {\bf Proof of Theorem 1.} We follow the loop: 
$(i)$ $\Rightarrow$ $(ii)$, $(ii)$ $\Rightarrow$ $(iii)$, $(iii)$ $\Rightarrow$ $(iv)$, $(iv)$ $\Rightarrow$ $(i)$. 
Then, we'll prove 
$(iii)$ $\Rightarrow$ 
s.R.e. is in $1$st class. 
\par
$(i)$ $\Rightarrow$ $(ii)$. Apply Lemma A.4 in $\S 9$ to the arithmetic function $F(h)=C_{f,g}(N,h)$. 
\par
$(ii)$ $\Rightarrow$ $(iii)$. Expand in finite Ramanujan expansion $g$, with a finite support of $g'$, hence of $\widehat{g}$, not depending on $h$ (likewise for $f,\widehat{g}$), inside $C_{f,g}(N,h)$ (so may exchange $h-$sum): 
$$
{1\over x}\sum_{h\le x}C_{f,g}(N,h)c_{\ell}(h)=\sum_q \widehat{g}(q)\sum_{n\le N}f(n){1\over x}\sum_{h\le x}c_q(n+h)c_{\ell}(h)
$$
\par
\noindent
whence, passing to the limit, normalizing (by Euler's function $\varphi(\ell)$, here) and applying $(ii)$, 
$$
\widehat{C_{f,g}}(N,\ell)={1\over {\varphi(\ell)}}\sum_q \widehat{g}(q)\sum_{n\le N}f(n)\lim_{x\to \infty}{1\over x}\sum_{h\le x}c_q(n+h)c_{\ell}(h), 
$$
\par
\noindent
in which, writing $\1_{\wp}\defineq 1$ if $\wp$ is true, $\defineq 0$ otherwise and $j\in \Z_{\ell}^*$ to abbreviate $j\le \ell$, $(j,\ell)=1$, we have  
$$
{1\over x}\sum_{h\le x}c_q(n+h)c_{\ell}(h)={1\over x}\sum_{r\in \Z_q^*}e^{2\pi i nr/q}\sum_{j\in \Z_{\ell}^*}\sum_{h\le x}e^{2\pi i (r/q-j/\ell)h}= 
$$
$$
=\1_{q=\ell}c_{\ell}(n)+O\left({1\over x}\sum_{r\in \Z_q^*}\sum_{j\in \Z_{\ell}^*}(1-\1_{q=\ell}\1_{r=j}){1\over {\left\Vert r/q-j/\ell\right\Vert}}\right) 
=\1_{q=\ell}c_{\ell}(n)+o(1), 
$$
\par
\noindent
letting: $x\to \infty$, proves the orthogonality relations (discussed in [M] with more details) 
$$
\lim_{x\to \infty}{1\over x}\sum_{h\le x}c_q(n+h)c_{\ell}(h)=\1_{q=\ell}c_{\ell}(n), 
$$
\par
\noindent
giving at once $(iii)$. 
\par
$(iii)$ $\Rightarrow$ $(iv)$. Observe that $\widehat{g}$ support is finite and independent of $h$. 
\par
$(iv)$ $\Rightarrow$ $(i)$. Trivial (since uniform convergence follows by $h-$independence). 
\par
\noindent
We now prove that $(iii)$ gives our s.R.e. is in $1$st class. 
\par
$(iii)$ $\Rightarrow$ applying bounds $(14)$, with $\varphi(\ell)\gg \ell/\log \ell$, see [CojM], 
$$
\left|\widehat{C_{f,g}}(N,\ell)\right|\le {{\left|\widehat{g}(\ell)\right|}\over {\varphi(\ell)}}\sum_{n\le N}\left|f(n)c_{\ell}(n)\right|
\ll_{\varepsilon} (\ell N)^{\varepsilon}{1\over {\ell^2}}\sum_{n\le N}(n,\ell), 
$$
\par
\noindent
from Lemma A.1 inequality $|c_{\ell}(n)|\le (n,\ell)$; since 
$$
\sum_{n\le N}(n,\ell)=\sum_{t|\ell}t\sum_{{n\le N}\atop {(n,\ell)=t}}1
\le \sum_{t|\ell}t\sum_{{n\le N}\atop {n\equiv 0\bmod t}}1
\le \sum_{t|\ell}t\left({N\over t}+1\right)
\le (N+\ell)d(\ell)
\ll_{\varepsilon} \ell^{\varepsilon}N, 
$$
\par
\noindent
from $\ell \ll N+h\ll N$ and the known bound on the divisor function $d(\ell)\ll_{\varepsilon} \ell^{\varepsilon}$, we get 
$$
\widehat{C_{f,g}}(N,\ell)\ll_{\varepsilon} {{N^{1+\varepsilon}}\over {\ell^2}}, 
$$
\par				% PAGE 11 
\noindent
i.e., by definition, the s.R.e. is in $1$st class. 
\par
\noindent
This completes the proof. \hfill $\square$ 

\bigskip
\bigskip

\par
\noindent {\bf Proof of Corollary 1.} We only need to prove the first formula, for which Lemma A.3 gives 
$$
\sum_{n\le N}f(n)c_{\ell}(n)=\widehat{f}(\ell)\varphi(\ell)N+O_{\varepsilon}\left((D\ell)^{1+\varepsilon}\right), 
$$
\par
\noindent
whence, from the explicit formula in Theorem 1 (recall $\ell,h,D\ll N$) again by Lemma A.1, 
$$
C_{f,g}(N,h)=\sum_{\ell \ll N}{{\widehat{g}(\ell)}\over {\varphi(\ell)}}\left(\widehat{f}(\ell)\varphi(\ell)N\right)c_{\ell}(h) 
 +O_{\varepsilon}\left(N^{\varepsilon}\sum_{\ell \ll N}{D\over {\ell}}(\ell,h)\right) 
$$
$$
=\sum_{\ell}\widehat{f}(\ell)\widehat{g}(\ell)c_{\ell}(h)N 
 +O_{\varepsilon}\Big(N^{\varepsilon}D\sum_{t|h}\sum_{v \ll N}{1\over v}\Big) 
=\SingSer_{f,g}(h)N+O_{\varepsilon}\Big(N^{\varepsilon}D\Big). 
$$
This completes the proof.\hfill $\square$ 

\bigskip
\bigskip

\par
\noindent
The shift-Ramanujan expansion, for any pair of arithmetic functions $f$ and $g$, leads us to a kind of entanglement of the two single Ramanujan expansions of, resp., $f$ and $g$ ! Our heuristic formul\ae \enspace (and others, in the literature), then, may be seen as a kind of squeezing on the diagonal, as obtained when we consider the same moduli $r=s$ in the single Ramanujan expansions, resp., with coefficients $\widehat{f}(r)$ and $\widehat{g}(s)$. This \lq \lq reduction on the diagonal\rq \rq, say, is a consequence, for our results [MS], [CMS], of the decay bounds for the single Ramanujan coefficients. However, as the above definition for the decay of this time the shift-Ramanujan coefficients reveals, this kind of reduction may hold in more general hypotheses (compare Corollary 1), than the ones we applied, now and in our previous studies (like, esp., [GMP],[MS],[CMS]). In particular,  in Theorem 1, the possibility to apply Carmichael's formula (implying $(iii)$, the explicit formula for the shift-Ramanujan expansion) seems to be the easiest requirement; we hope this will shed some light on the possibility to prove (in suitable, new hypotheses) the heuristic formul\ae, like $(3)$, for our shifted convolution sums. 

\medskip

\par
We wish, at this point, to conclude with three classical singular series, for $f,g$ correlations, thus giving (compare Corollary 1) heuristic formul\ae. 

\par
Of course, our first example is the case $f=g=\Lambda$ of $2k-$twin primes (there's a misprint in [CMS] at page 702) : 
$$%\label{2kTwins} (25)*
\SingSer_{\Lambda,\Lambda}(h)\defineq \sum_{q=1}^{\infty}{{\mu^2(q)}\over {\varphi^2(q)}}c_q(h), 
\leqno{(25)}
$$
\par
\noindent 
letting \enspace $h=2k$ \enspace to avoid $h$ odd (trivial case, with vanishing series). For it, we compare with the singular sum we get in case of truncations, say, $\Lambda_N$, i.e. 
$$%\label{finLambda} (26)*
\SingSer_{\Lambda_N,\Lambda_N}(h)\defineq \sum_{q=1}^{\infty}\widehat{\Lambda_N}(q)^2c_q(h) 
=\sum_{q\le N}\widehat{\Lambda_N}(q)^2c_q(h), 
\leqno{(26)}
$$
\par
\noindent 
that, see $(20)$ and $(21)$, has, in the range $q\le \sqrt{N}$, 
$$
\widehat{\Lambda_N}(q)={{\mu(q)}\over {\varphi(q)}}+O\left({1\over q}\exp\left(-c\sqrt{\log N}\right)\right)
\thinspace \Rightarrow \thinspace 
\widehat{\Lambda_N}^2(q)={{\mu^2(q)}\over {\varphi^2(q)}}+O\left({1\over {\varphi^2(q)}}\exp\left(-c\sqrt{\log N}\right)\right)
$$
\par
\noindent 
rendering 
$$
\SingSer_{\Lambda_N,\Lambda_N}(h)=\sum_{q\le \sqrt{N}}{{\mu^2(q)}\over {\varphi^2(q)}}c_q(h)+O\left(\exp\left(-c\sqrt{\log N}\right)\sum_{q\le \sqrt{N}}{{(q,h)}\over {\varphi^2(q)}}\right)+O\left(\sum_{q>\sqrt{N}}{{(q,h)}\over {\varphi^2(q)}}\right), 
$$
\par				% PAGE 12 
\noindent 
from Lemma 1; thus we may approximate (with a changed $c>0$) as 
$$%\label{SingSumLambda} (27)
\SingSer_{\Lambda_N,\Lambda_N}(h)=\SingSer_{\Lambda,\Lambda}(h)
 +O_h\left(\exp\left(-c\sqrt{\log N}\right)\right)+O_h\left({{\log^2 N}\over {\sqrt N}}\right). 
\leqno{(27)}
$$
\par
\noindent 
We use here the trivial bound, compare [CojM], $\varphi(q)\gg q/\log q$, inside the estimates 
$$
\sum_{q\le \sqrt{N}}{{(q,h)}\over {\varphi^2(q)}}\ll \sum_{d|h}d\sum_{{q\le \sqrt{N}}\atop {q\equiv 0\bmod d}}{{\log^2 q}\over {q^2}}
\ll (\log^2 N)\sum_{d|h}d^{-1}\sum_{n\le {{\sqrt{N}}\over d}}{1\over {n^2}}
\ll_h \log^2 N 
$$
\par
\noindent 
and 
$$
\sum_{q>\sqrt{N}}{{(q,h)}\over {\varphi^2(q)}}\ll \sum_{d|h}d\sum_{{q>\sqrt{N}}\atop {q\equiv 0\bmod d}}{{\log^2 q}\over {q^2}}
\ll \sum_{d|h}d^{-1}\sum_{n>{{\sqrt{N}}\over d}}{{\log^2 n+\log^2 d}\over {n^2}} 
$$
$$
\ll \sum_{d|h,d\le \sqrt{N}}{1\over d}\left(d{{\log^2 N}\over {\sqrt{N}}}+\sum_{n>N}{{\log^2 n}\over {n^2}}\right) 
    + \sum_{d|h,d>\sqrt{N}}{{\log^2 d}\over d}
\ll_h {{\log^2 N}\over {\sqrt N}}.
$$
\par
\noindent 
Thus $(27)$ proves that, at least in the range $q\le \sqrt N$ for present case, the \lq \lq singular sum\rq \rq \thinspace well approximates the original singular series. This, in case $f=g=\Lambda$ (for which we don't have absolute convergence of $\Lambda$ Ramanujan expansion, see [CMS]); but when we also have the absolute convergence of the original Ramanujan expansions, of both $f$ and $g$ separately, (so we are considering the single two Ramanujan expansions), we get an even better approximation (as the diligent reader may check previous two examples). 
\par
In fact, if we consider the other two examples (of [CMS] Corollaries) given in the previous $\S 4$, we see a very useful convergence of Ramanujan expansions.   We refer to our previous paper for the expression of the corresponding singular series (which are too involved to quote, due to space reasons). As the diligent reader may check, the singular series of these two examples (see $\S 4$ end) converge even better than our previous estimates, so the difference between them and finite Ramanujan expansions counterparts behaves much better than what we saw in $(27)$ (whose error terms depend on zero-free regions for the Riemann $\zeta$ function). 

\bigskip
\bigskip

\par
\centerline{\stampatello 6. Sieve functions}%\label{SieveF}6 
\bigskip
\par
\noindent
We first recall the definition [CL1,CL2,CL3]: a {\it sieve function} $f:\N \rightarrow \C$ of range $Q$ (that is an unbounded parameter, depending on other variables, see the following) may be written as 
$$
f(n)\defineq \sum_{d|n,d\le Q}f'(d), 
$$
\par
\noindent
where its Eratosthenes transform $f'$ (see the above) is {\it essentially bounded}, namely we recall 
$$
f'(q)\ll_{\varepsilon} q^{\varepsilon},
\qquad \hbox{\rm as} \enspace 
q\to \infty. 
$$
\par
\noindent
In other words, a sieve $f$ (of range $Q$) is a ($Q-$)truncated divisor sum satisfying the Ramanujan Conjecture. 
In fact, \lq \lq $f$ is essentially bounded\rq \rq, by M\"obius inversion, is equivalent to \lq \lq $f'$ is essentially bounded\rq \rq. 
\par
Hence, we immediately get the {\it finite Ramanujan expansion} for a sieve $f$ of range $Q$ 
$$
f(n)=\sum_{q\le Q}\widehat{f}(q)c_q(n), 
\quad \hbox{\rm where} \quad
\widehat{f}(q)=\sum_{d\equiv 0 \bmod q}{{f'(d)}\over d}, 
$$
\par
\noindent
with (at most) $Q$ terms. 
\par				% PAGE 13 
Sieve functions $f$ always have a mean-value (i.e., $M(f)\defineq \lim_x {1\over x}\sum_{n\le x}f(n)$, [M]) and it is $\widehat{f}(1)$ (this by Wintner's 1943 Theorem, [M], Theorem 2). Also, the Dirichlet series of a sieve $f$ of range $Q$ is the product of the Riemann $\zeta$ function and a Dirichlet polynomial with (at most) $Q$ terms : compare [CL3] (soon after $(1.3)$ equation). 
\par
From our considerations in $\S 2$, a general arithmetic function $f$ may be seen as a kind of sieve function of range $N$ (if we confine to its correlations). However, there is no point in considering the range $N$, as it holds for all arithmetic functions $f$. 
\par
\noindent
Note that the parameter $Q\ll N$  is dependent on $N$, or even on $h$ and, however, $Q\to \infty$, as $N\to \infty$ (avoiding trivialities). 
Hereafter, we assume $f'(Q)\neq 0$ and define the {\it level} \footnote{$^1$}{For a discussion on links between present definition and the {\it level of distribution} of $f$ {\it in arithmetic progressions}, compare [CL3].}
 of our $f$ as \enspace $\lambda(f)\defineq (\log Q)/(\log N)$. Notice the sensitivity of this definition to the main variable $N\to \infty$. 
\par
Expanding $f$ of range $D$ and $g$ of range $Q\ge D$ in finite Ramanujan expansions, we get 
$$%\label{RamaEntangle} (28)
C_{f,g}(N,h)=\sum_{d\le D}\widehat{f}(d)\sum_{q\le Q}\widehat{g}(q)\sum_{n\le N}c_d(n)c_q(n+h), 
\leqno{(28)}
$$
\par
\noindent
that we study by a modified large sieve inequality using the fact that 
Farey fractions ${j\over q}\neq {r\over d}$ are well-spaced in $[0,1]$:
$$
\sum_{n\le N}c_d(n)c_q(n+h)=\sum_{j\in \Z_q^*}e^{2\pi i hj/q}\sum_{r\in \Z_d^*}\sum_{n\le N}e^{-2\pi i nr/d}e^{2\pi i nj/q}= 
$$
$$
=\1_{d=q}Nc_q(h)+O\Big(\sum_{j\in \Z_q^*}\sum_{r\in \Z_d^*}\1_{{j\over q}\neq {r\over d}}{1\over {\left\Vert {j/q-r/d}\right\Vert}}\Big) 
$$
\par
\noindent
and the bounds on the Ramanujan coefficients, coming from $(13)$ in $\S 3$ and $f'$ and $g'$ bounds: 
$$%\label{sievefReBounds} (29)
\widehat{f}(d)\ll_{\varepsilon} D^{\varepsilon}/d
\quad \hbox{\rm and} \quad 
\widehat{g}(q)\ll_{\varepsilon} Q^{\varepsilon}/q
\leqno{(29)}
$$
\par
\noindent
together give (by Lemma 2 of [CSal] for our $DQ-$spaced Farey fractions $\lambda_r:=r/d$, $\lambda_s:=j/q$) 
$$
C_{f,g}(N,h)=N\sum_{q\le D}\widehat{f}(q)\widehat{g}(q)c_q(h) 
$$
$$
+O_{\varepsilon}\Big(D^{\varepsilon}\sum_{d\le D}{1\over {d^2}}\sum_{j\in \Z_d^*}\sum_{{r\in \Z_d^*}\atop {r\neq j}}{1\over {\left\Vert {{j-r}\over d}\right\Vert}}
  +Q^{\varepsilon}\sum_{d\le D}\sum_{{q\le Q}\atop {q\neq d}}{1\over {dq}}\sum_{j\in \Z_q^*}\sum_{r\in \Z_d^*}{1\over {\left\Vert {j/q-r/d}\right\Vert}}\Big)
$$
$$
=\SingSer_{f,g}(h)N 
 +O_{\varepsilon}\Big(D^{\varepsilon}\sum_{d\le D}{{\varphi(d)}\over {d^2}}\sum_{\Delta \le d/2}{d\over {\Delta}} 
  +DQ^{1+\varepsilon}\Big({\sum_{d\le D}{1\over {d^2}}\sum_{r\in \Z_d^*}1}\Big)^{1/2}\Big({\sum_{q\le Q}{1\over {q^2}}\sum_{j\in \Z_q^*}1}\Big)^{1/2}\Big),
$$
\par
\noindent
whence 
$$%\label{fRangeDgRangeQ} (30)
C_{f,g}(N,h)=\SingSer_{f,g}(h)N 
 +O_{\varepsilon}\left(D^{1+\varepsilon}+(DQ)^{1+\varepsilon}\right) 
=\SingSer_{f,g}(h)N 
 +O_{\varepsilon}\left((DQ)^{1+\varepsilon}\right). 
\leqno{(30)}
$$
\par
\noindent
This is an asymptotic formula, if $DQ\ll N^{1-\delta}$, for a fixed $\delta>\varepsilon$; i.e., calling $\lambda(f)$ and $\lambda(g)$, resp., the levels of, resp. $f$ and $g$, the requirement is $\lambda(f)+\lambda(g)<1$. In the particular case $f=g$ ($f$ autocorrelation) this means $\lambda(f)<1/2$ which is
 the well-known barrier for the large sieve technique (which we do not apply here but we rely actually on its proof as the main ingredient, namely the well-spaced property of Farey fractions). Notice the uniformity in $h\ge 0$ (uniformity with
respect to the shift which is useful for the correlation asymptotic formul\ae). 
\par				% PAGE 14 
Also, we may prove this formula by rendering rigorous the heuristic argument for fractional parts, given in $\S 2$, simply using 
the fact that they are
bounded; this gives an alternative to the present Farey fractions argument
and a much shorter proof. (The interested reader may fill in the details.) 

\bigskip
\bigskip

\par
\centerline{\stampatello 7. Sifting from small prime divisors}%\label{Sifted}7
\bigskip
\par
\noindent
We give a new definition which will be useful, when applying our study to sieve functions that, in some sense (we specify now), have no divisors with \lq \lq small primes\rq \rq. 
\par
We say that a general $f:\N \rightarrow \C$ satisfies a {\it sieve condition up to} $G$ ($>1$ and integer), by definition, when its Eratosthenes transform $f'$ (recall, defined by $f\ast \mu$) has the property 
$$%\label{S.c.G} (31)*
p\le G, \thinspace p|q 
\enspace \Longrightarrow \enspace 
f'(q)=0. 
\leqno{(31)}
$$
\par
\noindent
Notice that, in particular, since \thinspace $f(n)=\sum_{q|n}f'(q)$, defining (as usual, in sieve methods) the product of all primes up to $G$ as 
$$
P(G)\defineq \prod_{p\le G}p, 
$$
\par
\noindent
we obtain (for $f$ with a sieve condition up to $G$) a kind of sifting $f$ from small primes $p\le G$,  the formula 
$$
f(n)=f\left({n\over {(n,P(G))}}\right), 
\qquad \forall n\in \N. 
$$

\medskip

\par
We call $f$ a {\it $G-$sifted function} of range $Q$, whenever it is a sieve function of range $Q$ and satisfies a sieve condition up to $G$. 

\medskip

\par
Notice that from the formula $(13)$ for Ramanujan coefficients, we get that any $G-$sifted function $f$ (of range $Q$) has a finite Ramanujan expansion without indices from $2$ up to $G$ : 
$$
f(n)=\widehat{f}(1)+\sum_{G<q\le Q}\widehat{f}(q)c_q(n). 
$$
\par
\noindent
Also, the coefficients $\widehat{f}(q)=0$ when $q$ has a prime factor $p\le G$.  More precisely, 
$$%\label{fReGsifted} (32)*
f \enspace \hbox{\rm is} \thinspace \thinspace G-\hbox{\rm sifted} \enspace \hbox{\rm of \thinspace range } Q
\enspace \Longrightarrow \enspace 
f(n)=\widehat{f}(1)+\sum_{{G<q\le Q}\atop {(q,P(G))=1}}\widehat{f}(q)c_q(n). 
\leqno{(32)}
$$

\medskip

\par
\noindent
In particular, we have that the singular series $\SingSer_f(h)=\SingSer_{f,f}(h)$ (taking $g=f$ for the heuristic of $C_f=C_{f,f}$), in case $h>0$, of a $G-$sifted $f$ has the shape (using $(29)$ here) 
$$
\SingSer_f(h)=\widehat{f}(1)^2+\sum_{{G<q\le Q}\atop {(q,P(G))=1}}\widehat{f}(q)^2 c_q(h) 
$$
$$
=\widehat{f}(1)^2 + O_{\varepsilon}\left(Q^{\varepsilon}\sum_{G<q\le Q}{1\over {q^2}}|c_q(h)|\right) 
=\widehat{f}(1)^2 + O_{\varepsilon}\left(Q^{\varepsilon}\max_{G\le A\ll Q}{1\over {A^2}}\sum_{A<q\le 2A}|c_q(h)|\right) 
$$
$$
=\widehat{f}(1)^2 + O_{\varepsilon}\left(Q^{\varepsilon}d(h)/G\right)
=\widehat{f}(1)^2 + O_{\varepsilon}\left((hQ)^{\varepsilon}/G\right), 
$$
\par
\noindent
applying a dyadic argument based on the following bound, for all integers $0\le A<B$ (use $|c_q(h)|\le (q,h)$, Lemma A.1 $\S 9$, here): 
$$
\sum_{A<q\le B}|c_q(h)|\le \sum_{l|h}l\sum_{{A<q\le B}\atop {(q,h)=l}}1
\le \sum_{l|h}l\sum_{{A<q\le B}\atop {q\equiv 0 \bmod l}}1 
$$
$$				% PAGE 15 
=\sum_{l|h}l\left(\left[{B\over l}\right]-\left[{A\over l}\right]\right)
\le \sum_{l|h, \thinspace l\le B}l\left({B\over l}-{A\over l}+1\right)
\le 2\sum_{l|h}B
\ll Bd(h). 
$$

\medskip

\par
\noindent
Heuristically speaking, we have a kind of general philosophy, for sifted functions: \lq \lq low\rq \rq, Ramanujan coefficients vanish (like, on the other side, those out of range). The absence of low primes, $p\le G$, in conjunction with low shifts, again up to $G$, gives the following interesting properties, like $(33)$ \& $(34)$. 
\par
As an example, let us consider first what happens in arithmetic progressions, in case of sieve functions: Lemma A.2, $\S 9$, says that for a sieve $f$ of range $D\ll N$ we have 
$$
\sum_{{n\le N}\atop {n\equiv a\bmod t}}f(n)={N\over t}\sum_{k|t}\widehat{f}(k)c_k(a)+O_{\varepsilon}\left(D^{1+\varepsilon}\right); 
$$
\par
\noindent
so assume $t$ is $G-$sifted (a classical expression to mean $(t,P(G))=1$, here) and say the shift $h=-a$ satisfies $0<|a|=|h|\le G$: then, any prime divisor $p$ of $k|t$ in the right hand side above has to be $p>G$, while any prime divisor of $a$ cannot be greater than $G$ itself.  Then the previous formula simplifies to (using $\widehat{f}(k)\ll_{\varepsilon} D^{\varepsilon}/k$ and the fact that $k|t$ and $k>1$ $\Rightarrow$ $k>G$)
$$
\sum_{{n\le N}\atop {n\equiv a\bmod t}}f(n)={N\over t}\sum_{k|t}\widehat{f}(k)\mu(k)+O_{\varepsilon}\left(D^{1+\varepsilon}\right)
={N\over t}\widehat{f}(1)+O_{\varepsilon}\left(D^{\varepsilon}\left({{Nt^{\varepsilon}}\over {tG}}+D\right)\right), 
$$
\par
\noindent
entailing 
$$%\label{sumAPh} (33)
t \enspace \hbox{\rm is} \enspace G-\hbox{\rm sifted}, \enspace 0<|h|\le G 
\enspace \Longrightarrow \enspace 
\sum_{{n\le N}\atop {n\equiv -h\bmod t}}f(n)={N\over t}\widehat{f}(1)+O_{\varepsilon}\left(D^{\varepsilon}\left({{Nt^{\varepsilon}}\over {tG}}+D\right)\right), 
\leqno{(33)}
$$
\par
\noindent
which is, for what we specified, uniform in the non-zero residue classes $h\in [-G,G]$. 
\par
These same hypotheses for $a=-h$ give furthermore, from Lemma A.3 of $\S 9$, 
$$
\sum_{n\le N}f(n)c_{\ell}(n-a)=\mu(\ell)\widehat{f}(\ell)N+O_{\varepsilon}\left((D\ell)^{1+\varepsilon}\right), 
$$
\par
\noindent
uniformly in $0<|a|\le G$, whenever, this time, $\ell$ is $G-$sifted; as before, for $f$ of range $D$, 
$$%\label{sumclh} (34)
\ell \enspace \hbox{\rm is} \enspace G-\hbox{\rm sifted}, \enspace 0<|h|\le G 
\enspace \Longrightarrow \enspace 
\sum_{n\le N}f(n)c_{\ell}(n+h)=\mu(\ell)\widehat{f}(\ell)N+O_{\varepsilon}\left((D\ell)^{1+\varepsilon}\right), 
\leqno{(34)}
$$
\par
\noindent
uniformly in shifts $0<|h|\le G$. 

These formul\ae \enspace are very useful inside the correlations. In fact, these may be expressed as  
$$
C_{f,g}(N,h)=\sum_{n\le N}f(n)g(n+h)=\sum_{q\le N+h}g'(q)\sum_{{n\le N}\atop {n\equiv -h\bmod q}}f(n), 
$$
\par
\noindent
in which, if $f$ has range $D\ll N$ and $g$ is $G-$sifted of range $Q\ll N$, all $q$ are $G-$sifted, so 
$$
C_{f,g}(N,h)=\sum_{n\le N}f(n)g(n+h)
=\sum_{q\le Q}g'(q)\widehat{f}(1){N\over q}
 +O_{\varepsilon}\left(N^{\varepsilon}\sum_{q\le Q}|g'(q)|\left(D+{N\over {qG}}\right)\right), 
$$
\par
\noindent
from $(33)$, for the $0<h\le G$, therefore with $G=o(N)$ we have 
$$
C_{f,g}(N,h)=\sum_{n\le N}f(n)g(n+h)
=\widehat{f}(1)\widehat{g}(1)N+O_{\varepsilon}\left(N^{\varepsilon}\left(DQ+N/G\right)\right). 
$$
\par				% PAGE 16 
\noindent
Compared to the results $(30)$ in $\S 6$, apart from a short shift (and sifting hypothesis), we have now in the main term a collapse to the lonely $\widehat{f}(1)\widehat{g}(1)$ (the first of $\SingSer_{f,g}(h)$, since $c_1(h)=1$). 
\par
This is in accordance with heuristics in case $f=g$ above: in fact, our singular series is now (from the hypothesis $g$ is $G-$sifted) 
$$
\SingSer_{f,g}(h)=\widehat{f}(1)\widehat{g}(1)+\sum_{G<\ell \le Q}\widehat{f}(\ell)\widehat{g}(\ell)
=\widehat{f}(1)\widehat{g}(1)+O_{\varepsilon}\left(N^{\varepsilon}/G\right). 
$$

\par
\noindent
By the way, we might have used, also, the other formula, $(34)$, to estimate our correlation: 
$$
C_{f,g}(N,h)=\sum_{q\le Q}\widehat{g}(q)\sum_{n\le N}f(n)c_q(n+h)
=\sum_{q\le Q}\widehat{g}(q)\mu(q)\widehat{f}(q)N 
 + O_{\varepsilon}\left(N^{\varepsilon}D\sum_{q\le Q}|\widehat{g}(q)|q\right) 
$$
$$
=\widehat{f}(1)\widehat{g}(1)N+\sum_{G<q\le Q}\widehat{g}(q)\mu(q)\widehat{f}(q)N 
 + O_{\varepsilon}\left(N^{\varepsilon}DQ\right)
=\widehat{f}(1)\widehat{g}(1)N
 + O_{\varepsilon}\left(N^{\varepsilon}(DQ+N/G)\right), 
$$
\par
\noindent
where we are still using hypotheses $0<|h|\le G$, $G=o(N)$, $g$ is $G-$sifted of range $Q\ll N$ and $f$ sieve of range $D\ll N$. In this case we have again the remainder $\ll_{\varepsilon}N^{1+\varepsilon}/G$, that \lq \lq can't be removed\rq \rq. 
\par
Anyway, this parameter $G\to \infty$, $G=o(N)$ as $N\to \infty$, in our formul\ae, has to be chosen (at least) as $G\gg N^{\delta}$, for some fixed $\delta>0$. (Compare the Remark after Lemma A.5 in $\S 9$.)

\bigskip

\par
We refer the interested reader to the last section $\S 9$, for further properties, involving sieve functions and numbers free of small prime factors, say without prime divisors $p\le G$ : especially, to Lemma A.5. 
\par
In fact, this Lemma shows that, in sums of sieve functions (with a shift $0<|h|\le G$), the condition of being coprime to a fixed natural number $q$,  which is free of prime factors $p\leq G$, may be dropped at a small cost (depending, of course, on $G$). 
\par
We leave for our reader, a simple exercise about this property: with the above considerations starting from $(33)$, apply third formula of Lemma A.5, $\S 9$, in order to prove the noteworthy property, for the, say, $q-${\it coprime correlation} of $f$ and $g$ (for which recover above hypotheses), namely (compare both above results for usual correlation) 
$$
C_{f,g}^{(q)}(N,h)\defineq \sum_{{n\le N}\atop {(n,q)=1}}f(n)g(n+h)
=\widehat{f}(1)\widehat{g}(1)N
 + O_{\varepsilon}\left(N^{\varepsilon}(DQ+N/G)\right). 
$$
\par
\noindent
This time the error term \enspace $N^{1+\varepsilon}/G$ \enspace is also the cost we pay to \lq \lq remove\rq \rq   the condition of coprimality with $q$, which is free of factors $p\le G$. 

\par

(We see that, when considering averages of correlations, esp. [CL1], the really important part of our singular series is the first term, like here with the sifting hypothesis. In some sense, averaging over the shift, comparing the symmetry integral calculations in $\S 9$, \lq \lq smooths\rq \rq, so to speak, the arithmetic irregularities; which are \lq \lq overridden\rq \rq, here, by the \lq \lq no small primes condition\rq \rq.)

\bigskip
\bigskip

\par
\centerline{\stampatello 8. Concluding Remarks}%\label{CoRemarx}8
\bigskip
\par
\noindent
In our previous paper [CMS2] we introduced the main idea entailing finite Ramanujan expansions, inside the shifted convolution sums of arbitrary arithmetic functions $f$ and  $g$ (see $(12)$ in section 3): namely, turning $f$ and $g$ into their truncated divisor sums counterparts, see $(11)$ in section 2. 
\par
This is of particular importance, both because finite Ramanujan expansions (apart from Hildebrand's quoted examples, which are not pure in our sense) have never been studied before and, on the other hand, the truncated divisor sums are very well known, especially in sieve theory (compare the \lq \lq sieve functions\rq \rq, in section 6). The present paper gives some properties and examples for f.R.e., in sections 3 and 4; actually, we are only starting
here a general theory of general f.R.e. 
\par				% PAGE 17 
The present paper's new idea is another kind of Ramanujan expansion (not regarding general single arithmetic functions), the one for the shifted convolution sum of $f$ and $g$ (as an arithmetic function, with argument the shift itself), called the \lq \lq shift-Ramanujan expansion\rq \rq. 
This is a very special kind of Ramanujan expansion, \lq \lq entangling\rq \rq \thinspace the two single Ramanujan expansions (we might say: marries them!), of $f$ \& $g$ (like, by definition, the shifted convolution sum of $f$ and $g$ entangles them). This is clear in the proof of our Theorem, given in section 5, and is evident in the formul\ae, for the shifted convolution sum of $f$ \& $g$ in terms of the (finite-)Ramanujan coefficients $\hat{f}$ \& $\hat{g}$ of, resp., $f$ \& $g$ : see $(28)$. 
\smallskip
\par
In fact, while in section 6 these formul\ae \enspace give a duality  between $f$ and $g$ entanglement (in their correlation) and $\widehat{f}$ and $\widehat{g}$ entanglement (see formula $(28)$, esp.), we already see a duality in between $f$ and $g$ correlation and $f'$ and $g'$  entanglement (with the introduction of fractional parts, outside main term) in section 2 (see where we start giving heuristic formul\ae \enspace for the $f$ and $g$ correlation). 
\smallskip
\par
By the way, the word \lq \lq duality\rq \rq, here, used freely in our context, has a resemblance in links  between different spectral analyses for a problem; in fact, the \lq \lq spectral formula\rq \rq, for correlation, needed to prove the heuristic in section 2, is based exactly on the (finite) Fourier expansion for the fractional parts, coming, as usual, from $1-$periodic Bernoulli function. Of course, we also have two kinds of Ramanujan expansions (say, two \lq \lq harmonic formul\ae\rq \rq) : the one for the single $f$ and  $g$, with coefficients $\hat{f}$ and $\hat{g}$, featuring correlations of Ramanujan sums (with 2 moduli, i.e. $(28)$ in section 6) and, in parallel, the shift-Ramanujan expansion with a kind of \lq \lq mysterious\rq \rq, new, shift-Ramanujan coefficients times the Ramanujan sum (with 1 modulus). 
\smallskip
\par
These two ways of expanding the correlation of $f$ \& $g$ give two possible approaches: the one with (single \&) finite Ramanujan expansions entangle the two moduli (of $\hat{f}$ \& $\hat{g}$) inside the correlation of Ramanujan sums w.r.t. the same two moduli; while, the shift-Ramanujan expansion entangles the two functions $f$ \& $g$ (so, the moduli of $\hat{f}$ \& $\hat{g}$) inside the \lq \lq black box\rq \rq, given by the shift-Ramanujan coefficients, $\widehat{C_{f,g}}$ (so, at last, shift-Ramanujan expansions have only one modulus, apparently). In some sense, our Theorem tries to take a glance (in suitable hypotheses) to this black box ! 

\medskip

\par
The \lq \lq spectral analysis\rq \rq, we are dealing with, of course, is elementary here, nothing sophisticated like the one for shifted convolution sums, say, in the context of modular forms as in the Rankin-Selberg method. (Even if we think that there are, almost surely, links to that: compare, for example, our formula for the shifted convolution sums, with a coprimality condition, we have in section 7 end.) 

\bigskip
\bigskip
\bigskip

\par
\centerline{\stampatello 9. Appendix}%\label{App}9
\bigskip
\par
\noindent
The results here are listed in increasing order of technicality. 
\bigskip
\par
We start, as we use this bound for the next Lemma, with the following elementary property (which we could not locate in the literature). 
\smallskip
\par
\noindent {\bf Lemma A.0.} {\it For all } $a,b\in \N$ {\it we have} 
$$
\varphi(ab)\le a\varphi(b). 
$$
\medskip
\par
\noindent {\bf Remark.} A straightforward proof comes considering lattice points of the $a$ times $b$ rectangle. 
\medskip
\par
\noindent {\bf Proof.} We use the well-known 
$$
\varphi(n)=n\prod_{p|n}\left(1-{1\over p}\right) 
\quad 
\forall n\in \N, 
$$
\par
\noindent
in order to get 
$$
{{\varphi(ab)}\over {\varphi(b)}}=a\prod_{p|ab}\left(1-{1\over p}\right)\prod_{p|b}\left(1-{1\over p}\right)^{-1} 
=a\prod_{{p|ab}\atop {b\not \equiv 0\bmod p}}\left(1-{1\over p}\right), 
$$
\par
\noindent
whence we have the lemma, since the last product is $1$ if empty, otherwise $<1$ (like all factors).\enspace $\square$ 

\bigskip
\bigskip

\par				% PAGE 18 
We prove, as we often use this bound in the above sections, the following Lemma. 
\smallskip
\par
\noindent {\bf Lemma A.1.} {\it For all } $q\in \N$ {\it and } $n\in \Z$ {\it we have} 
$$
|c_q(n)|\le (n,q). 
$$
\smallskip
\par
\noindent {\bf Proof.} We simply put together H\"older's identity [M] $c_q(n)=\varphi(q)\mu(q/(n,q))/\varphi(q/(n,q))$ and previous Lemma, with $a=(n,q)$, $b=q/(n,q)$, to get the inequality 
$$
\left|c_q(n)\right|=\left|{{\varphi(q)\mu(q/(n,q))}\over {\varphi(q/(n,q))}}\right|
\le {{\varphi(q)}\over {\varphi(q/(n,q))}}
\le (n,q).\enspace \square 
$$

\bigskip

\par
We state and prove a result that we need, but is also of independent interest. 
\smallskip
\par
\noindent {\bf Lemma A.2.} {\it Let } $D,N,t\in \N$ {\it and assume } $f:\N \rightarrow \C$ {\it is a sieve function of range } $D\ll N$. {\it Then, uniformly in } $a\in \Z$, 
$$
\sum_{{n\le N}\atop {n\equiv a\bmod t}}f(n)={N\over t}\sum_{k|t}\widehat{f}(k)c_k(a)+O_{\varepsilon}\left(D^{1+\varepsilon}\right). 
$$
\smallskip
\par
\noindent {\bf Remark.} We may clearly assume that $0\le |a|\le t$, whence $a\ll N$. 
\smallskip
\par
\noindent {\bf Proof.} Opening $f(n)$, our left hand side is 
$$
\sum_{d\le D}f'(d)\sum_{{n\le N}\atop {{n\equiv 0\bmod d}\atop {n\equiv a\bmod t}}}1
=\sum_{{d\le D}\atop {(d,t)|a}}f'(d)\left({{N(d,t)}\over {dt}}+O(1)\right), 
$$
\par
\noindent
where the remainder is 
$$
\ll_{\varepsilon} D^{\varepsilon}\sum_{{d\le D}\atop {(d,t)|a}}1
\ll_{\varepsilon} D^{1+\varepsilon}. 
$$
\par
\noindent
The main term here is 
$$
{N\over t}\sum_{{d\le D}\atop {(d,t)|a}}{{f'(d)(d,t)}\over d} 
$$
\par
\noindent
for which we evaluate (using $(1)$ and definitions above) 
$$
\sum_{{d\le D}\atop {(d,t)|a}}{{f'(d)(d,t)}\over d}=\sum_{{\ell|t}\atop {\ell|a}}\ell \sum_{{d_0}\atop {(d_0,t/\ell)=1}}{{f'(\ell d_0)}\over {\ell d_0}}
=\sum_{{\ell|t}\atop {\ell|a}}\ell\sum_{m\left|{t\over {\ell}}\right.}\mu(m)\sum_{d_1}{{f'(\ell md_1)}\over {\ell md_1}} 
$$
$$
=\sum_{{\ell|t}\atop {\ell|a}}\ell\sum_{m\left|{t\over {\ell}}\right.}\mu(m)\widehat{f}(\ell m)
=\sum_{k|t}\widehat{f}(k)\sum_{{\ell|k}\atop {\ell|a}}\ell \mu\left({k\over {\ell}}\right)
=\sum_{k|t}\widehat{f}(k)c_k(a).\enspace \square 
$$

\bigskip

\par
The previous lemma proves the following result, that has an independent interest, too. 
\smallskip
\par
\noindent {\bf Lemma A.3.} {\it Let } $\ell,D,N\in \N$ {\it and assume } $f:\N \rightarrow \C$ {\it is a sieve function of range } $D\ll N$. {\it Then, uniformly in } $a\in \Z$, 
$$
\sum_{n\le N}f(n)c_{\ell}(n-a)=\widehat{f}(\ell)c_{\ell}(a)N+O_{\varepsilon}\left((D\ell)^{1+\varepsilon}\right). 
$$
\smallskip
\par
\noindent {\bf Proof.} Inserting
$$
c_{\ell}(n-a)=\sum_{{t|\ell}\atop {t|n-a}}t\mu\left({{\ell}\over t}\right), 
$$
\par				% PAGE 19 
\noindent
(see $(1)$), the left hand side is (from the previous lemma) 
$$
\sum_{t|\ell}t\mu\left({{\ell}\over t}\right)\sum_{{n\le N}\atop {n\equiv a\bmod t}}f(n)
=N\sum_{t|\ell}\mu\left({{\ell}\over t}\right)\sum_{k|t}\widehat{f}(k)c_k(a)+O_{\varepsilon}\left((D\ell)^{1+\varepsilon}\right), 
$$
\par
\noindent
whence we simply get the result, from 
$$
\sum_{t|\ell}\mu\left({{\ell}\over t}\right)\sum_{k|t}\widehat{f}(k)c_k(a)
=\sum_{k|\ell}\widehat{f}(k)c_k(a)\sum_{t'\left|{{\ell}\over k}\right.}\mu\left({{\ell}\over {kt'}}\right)
=\widehat{f}(\ell)c_{\ell}(a).\enspace \square 
$$

\bigskip

\par
\noindent {\bf Remark.} We explicitly point out that it is very important to have any improvement in the remainders of these two lemmas (and it is clear that once A.2 has a better error bound, then A.3 also has). In fact, this is true for both the lemmas and for their applications, especially to the main  terms in explicit formul\ae, for correlations; in particular, the explicit formula for shift Ramanujan expansion coefficients (compare Theorem 1, section $1$) requires a straightforward application (in the proof of Corollary 1, see $\S 5$) of Lemma A.3. 

\bigskip

\par
We give a very useful lemma, esp., for the shift Ramanujan expansions. 
\smallskip
\par
\noindent {\bf Lemma A.4.} {\it Let } $F:\N\rightarrow \C$ {\it have an uniformly convergent Ramanujan expansion, i.e.} 
$$
F(h)=\sum_{q=1}^{\infty}\widehat{F}(q)c_q(h), 
\enspace 
\forall h\in \N, 
$$
\par
\noindent
{\it with some coefficients } $\widehat{F}(q)\in \C$ {\it independent of } $h$ ({\it even in their support}). {\it Then, these are} 
$$%\label{Carmichael} (35)
\widehat{F}(\ell)={1\over {\varphi(\ell)}}\lim_{x\to \infty}{1\over x}\sum_{h\le x}F(h)c_{\ell}(h). 
\leqno{(35)}
$$
\medskip
\par
\noindent {\bf Remark.} We call the above \lq \lq Carmichael's formula\rq \rq, for Ramanujan coefficients [Car], [M]. 
\medskip
\par
\noindent {\bf Proof.} 
Fix $\ell\in \N$ and, by uniform convergence, we have \enspace $\forall \varepsilon>0$ $\exists Q=Q(\varepsilon,\ell)$, with $Q>\ell$ and 
$$
\left|\sum_{q>Q}\widehat{F}(q)c_q(h)\right|<{{\varepsilon}\over {d(\ell)}}, 
$$
\par
\noindent
entailing 
$$
{1\over x}\sum_{h\le x}F(h)c_{\ell}(h)=\sum_{q\le Q}\widehat{F}(q){1\over x}\sum_{h\le x}c_{\ell}(h)c_q(h)+{1\over x}\sum_{h\le x}c_{\ell}(h)\sum_{q>Q}\widehat{F}(q)c_q(h) 
$$
\par
\noindent
(notice purity allows sums exchange) implies (\lq \lq ${\displaystyle \lim_x}$\rq \rq, here, abbreviating \lq \lq ${\displaystyle \lim_{x\to \infty}}$\rq \rq) 
$$
\left|{1\over {\varphi(\ell)}}\lim_x {1\over x}\sum_{h\le x}F(h)c_{\ell}(h)
 -{1\over {\varphi(\ell)}}\sum_{q\le Q}\widehat{F}(q)\lim_x {1\over x}\sum_{h\le x}c_{\ell}(h)c_q(h)\right|
  \le {{\varepsilon}\over {\varphi(\ell)d(\ell)}} \lim_x {1\over x}\sum_{h\le x}(\ell,h), 
$$
\par
\noindent
from \enspace $|c_{\ell}(h)|\le (\ell,h)$, see Lemma A.1, whence the orthogonality relations (see [M]) 
$$
\lim_x {1\over x}\sum_{h\le x}c_{\ell}(h)c_q(h)=\sum_{j\in \Z_{\ell}^*}\sum_{r\in \Z_q^*}\lim_x {1\over x}\sum_{h\le x}e^{2\pi i (j/\ell-r/q)h} 
=\sum_{j\in \Z_{\ell}^*}\sum_{r\in \Z_q^*}\1_{q=\ell}\1_{r=j}
=\1_{q=\ell}\varphi(\ell)
$$
\par				% PAGE 20 
\noindent
and the formula 
$$
{1\over x}\sum_{h\le x}(\ell,h)=\sum_{t|\ell}{t\over x}\sum_{{h'\le {x\over t}}\atop {(h',{{\ell}\over t})=1}}1
=\sum_{t|\ell}{t\over x}\sum_{d\left|{{\ell}\over t}\right.}\mu(d)\left[{x\over {dt}}\right]
=\sum_{t|\ell}\sum_{d\left|{{\ell}\over t}\right.}{{\mu(d)}\over d}+O\left({1\over x}\sum_{t|\ell}td(\ell/t)\right) 
$$
$$
=\sum_{t|\ell}{{\varphi(\ell/t)}\over {\ell/t}}+o(1) 
=\sum_{d|\ell}{{\varphi(d)}\over d}+o(1) 
$$
\par
\noindent
(flipping to the complementary divisor $d=\ell/t$), as $x\to \infty$, together give 
$$
\left|{1\over {\varphi(\ell)}}\lim_{x\to \infty}{1\over x}\sum_{h\le x}F(h)c_{\ell}(h)-\widehat{F}(\ell)\right|
 \le {{\varepsilon}\over {\varphi(\ell)d(\ell)}} \sum_{d|\ell}{{\varphi(d)}\over d}
 \le {{\varepsilon}\over {\varphi(\ell)}}
 \le \varepsilon, 
$$
\par
\noindent
which, as \enspace $\varepsilon>0$ \enspace is arbitrary, shows $(35)$.\enspace $\square$ 

\bigskip

\par
We give, now, an important result, with applications to the $G-$sifted functions. 
\smallskip
\par
\noindent {\bf Lemma A.5.} {\it Let } $G,q,N\in \N$ {\it and assume } $f:\N \rightarrow \C$ {\it is essentially bounded, while } $q$ {\it has all prime factors greater than } $G$, {\it i.e.}, \enspace $(q,P(G))=1$. {\it Then} 
$$
\sum_{{n\le N}\atop {(n,q)=1}}f(n)=\sum_{n\le N}f(n)+O_{\varepsilon}\left((qN)^{\varepsilon}{N\over G}\right). 
$$
\par
\noindent
{\it Furthermore, assuming } $0<|a|\le G$, 
$$
\sum_{{n\le N}\atop {{n\equiv a\bmod t}\atop {(n,q)=1}}}f(n)=\sum_{{n\le N}\atop {n\equiv a\bmod t}}f(n)
 +O_{\varepsilon}\left((qN)^{\varepsilon}\left({N\over {tG}}+1\right)\right), 
$$
\par
\noindent
{\it whence, adding the hypothesis } $f$ {\it is a sieve function of range $D\ll N$, from Lemma A.2,} 
$$
\sum_{{n\le N}\atop {{n\equiv a\bmod t}\atop {(n,q)=1}}}f(n)={N\over t}\sum_{k|t}\widehat{f}(k)c_k(a)
 +O_{\varepsilon}\left((qN)^{\varepsilon}\left({N\over {tG}}+D\right)\right); 
$$
\par
\noindent
{\it also, sieve $f$ of range $D\ll N$ have, in case $0<|a|\le G$,} 
$$
\sum_{{n\le N}\atop {(n,q)=1}}f(n)c_{\ell}(n-a)=\widehat{f}(\ell)c_{\ell}(a)N
 +O_{\varepsilon}\left((\ell qN)^{\varepsilon}\left({N\over G}+D\ell\right)\right) 
$$
\par
\noindent
{\it whence, applying our Lemma A.3 for case $0<|a|\le G$,} 
$$
\sum_{{n\le N}\atop {(n,q)=1}}f(n)c_{\ell}(n-a)=\sum_{n\le N}f(n)c_{\ell}(n-a)
 +O_{\varepsilon}\left((\ell qN)^{\varepsilon}\left({N\over G}+D\ell\right)\right). 
$$
\smallskip
\par
\noindent {\bf Proof.} Lemma 3 of [CMS2], $\1_{(n,q)=1}=\sum_{d|q,d|n}\mu(d)$, entails first LHS is (after, use: \hfil $d(q)\ll_{\varepsilon} q^{\varepsilon}$) 
$$
\sum_{d|q}\mu(d)\sum_{m\le {N\over d}}f(dm)=\sum_{n\le N}f(n)+O_{\varepsilon}\Big(N^{\varepsilon}\sum_{{d|q}\atop {d>G}}{N\over d}\Big)
=\sum_{n\le N}f(n)+O_{\varepsilon}\Big(N^{1+\varepsilon}{{d(q)}\over G}\Big).
$$
\par				% PAGE 21 
\noindent
In the same way 
$$
\sum_{{n\le N}\atop {{n\equiv a\bmod t}\atop {(n,q)=1}}}f(n)
=\sum_{d|q}\mu(d)\sum_{{n\le N}\atop {{n\equiv a\bmod t}\atop {n\equiv 0\bmod d}}}f(n) 
=\sum_{{n\le N}\atop {n\equiv a\bmod t}}f(n)+O_{\varepsilon}\Big(N^{\varepsilon}
                                                             \sum_{{d|q}\atop {{d>G}\atop {(d,t)|a}}}\Big({{N(d,t)}\over {dt}}+1\Big)
                                                              \Big), 
$$
\par
\noindent
but, now, $(d,t)=1$, since $(d,t)>1$ implies the primes $p|(d,t)$ are, dividing $a\neq 0$, all not greater than $G$, which is absurd by our assumption. Thus we have proved our second formula. In turn, this can be used to prove our fourth formula (compare Lemma A.3 proof), 
$$
\sum_{{n\le N}\atop {(n,q)=1}}f(n)c_{\ell}(n-a)=\sum_{t|\ell}t\mu\left({{\ell}\over t}\right)\sum_{{n\le N}\atop {n\equiv a\bmod t}}f(n)
 +O_{\varepsilon}\left((\ell qN)^{\varepsilon}(N/G+D\ell)\right) 
$$
$$
=\widehat{f}(\ell)c_{\ell}(a)N+O_{\varepsilon}\left((\ell qN)^{\varepsilon}(N/G+D\ell)\right).\enspace \square 
$$

\bigskip

\par
\noindent {\bf Remark.} Since $\varepsilon>0$ is arbitrarily small, $G\gg N^{\delta}$ ($\delta>0$ fixed) is preferred in all formul\ae. 

\medskip

\par
We express, here, the singular series in a very useful way, for heuristic calculations above. 
\smallskip
\par
\noindent {\bf Lemma A.6.} {\it Let } $h\ge 0$ {\it be integer and assume } $f,g:\N \rightarrow \C$ {\it have finite Ramanujan expansions. Then} 
$$
\SingSer_{f,g}(h)=\sum_{l|h}l\sum_{d}{{f'(d)}\over d}\sum_{{q}\atop {(q,d)=l}}{{g'(q)}\over q}. 
$$
\smallskip
\par
\noindent {\bf Proof.} In fact, $\SingSer_{f,g}$ definition, then $(1)$ and $(13)$ give (when $h=0$, $l|h$ means any $l\ge 1$): 
$$
\SingSer_{f,g}(h)=\sum_{q=1}^{\infty}\widehat{f}(q)\widehat{g}(q)c_q(h)=\sum_{q=1}^{\infty}\widehat{f}(q)\widehat{g}(q)\sum_{{l|h}\atop {l|q}}l \mu\left({q\over l}\right) 
=\sum_{l|h}l \sum_{k}\mu(k)\widehat{f}(lk)\widehat{g}(lk) 
$$
$$
=\sum_{l|h}l\sum_{k}\mu(k)\sum_{d}{{f'(lkd)}\over {lkd}}\sum_{q}{{g'(lkq)}\over {lkq}}
=\sum_{l|h}\doublesum_{(t,r)=1}{{f'(lt)g'(lr)}\over {ltr}}
$$
$$
=\sum_{l|h}l\doublesum_{(d,q)=l}{{f'(d)g'(q)}\over {dq}}.\enspace \square 
$$

\bigskip

\par
\noindent
We explore, now, the properties of a particular arithmetic function, so to clarify the properties of regularity (see soon after Theorem 1, in $\S 1$) for the shift-Ramanujan expansions. 

\medskip

\par
We build here an example of an arithmetic function, for which we prove that its s.R.e. is not in first class. It is denoted $f_H$, as it depends on the \lq \lq length\rq \rq, of short intervals $[x-H,x+H]$ for which we check the symmetry (in particular, for almost all of them, namely all $x\in (N,2N]$, safe at most $o(N)$ of them, as $N\to \infty$); we define it, generically, as a function assuming only two different values $c_1,c_2\in \C$, in consecutive intervals (of integers) with length $H$ (say, it starts $f_H(n)=c_1$, when $n\le H$, then $f_H(n)=c_2$, in $(H,2H]$, for example) and it is periodic of period $2H$; it clearly has mean-value $(c_1+c_2)/2$ and is not constant (as $c_1\neq c_2$). However, its \lq \lq short interval mean-value\rq \rq, say, has sometimes large size exactly when the short interval length is $H$; so, we don't consider its \lq \lq Selberg integral\rq \rq: 
$$
J_{f_H}(N,H)\defineq \sum_{N<x\le 2N}\Big|\sum_{x<n\le x+H}f_H(n)-M_{f_H}(x,H)\Big|^2, 
$$
\par
\noindent
as the difficulties show up already for the definition of this \lq \lq expected mean-value\rq \rq, in short intervals, $M_{f_H}(x,H)$. (Compare [CL1],[CL2], for the definitions and following properties.)
\par				% PAGE 22 
Thus we consider the \lq \lq easier\rq \rq \enspace symmetry integral ($\sgn(0)\defineq 0$, $r\neq 0\Rightarrow \sgn(r)\defineq {{|r|}\over r}$) 
$$
J_{f_H\, ,\,\sgn}(N,H)\defineq \sum_{N<x\le 2N}\Big|\sum_{x-H\le n\le x+H}\sgn(n-x)f_H(n)\Big|^2, 
$$
\par
\noindent
that checks the symmetry, around $x$, in {\it almost all} the short intervals $[x-H,x+H]$; namely, neglecting at most $o(N)$ of them for $x\in (N,2N]$ 
and the length is {\it short}, i.e. $H=o(N)$, as $N\to \infty$. The short intervals expected mean-value vanishes, here in the $f$ symmetry integral, whatever the $f:\N \rightarrow \C$ is (and this by definition). 
\par
We know (from Lemma 7 of [CL1]) that, since $f_H$ is bounded, hereafter $H=o(N)$, 
$$
J_{f_H\, ,\,\sgn}(N,H)=\sum_{h}W_H(h)\sum_{N<n\le 2N}f_H(n)f_H(n-h) + O(H^3), 
$$
\par
\noindent
where the weight $W_H(h)$ is the correlation of $\sgn$ function, [CL2], namely 
$$
W_H(h)\defineq \doublesum_{{-H\le h_1\thinspace ,\thinspace h_2\le H}\atop {h_2-h_1=h}}\sgn(h_1)\sgn(h_2), 
$$
\par
\noindent
which we will not write explicitly (see [CL2]), but only use its properties [CL2], like : $W_H$ is even, 
$$
W_H(h)\ll H,
\quad 
W_H(h)=0, \forall h\not \in [-2H,2H]. 
$$
\par
\noindent
We have a difference with our correlations, namely for $h<0$, abbreviating $C_{f_H}$ for $C_{f_H,f_H}$, 
$$
\sum_{N<n\le 2N}f_H(n)f_H(n-h)=C_{f_H}(2N,-h)-C_{f_H}(N,-h), 
$$
\par
\noindent
while the case $h>0$ (and $h\le 2H$, here) is handled from the immediate change of variables 
$$
\sum_{N<n\le 2N}f_H(n)f_H(n-h)=\sum_{N-h<n\le 2N+h}f_H(n)f_H(n+h) 
=C_{f_H}(2N,h)-C_{f_H}(N,h)+O(H) 
$$
\par
\noindent
and the case $h=0$ is negligible: 
$$
\sum_{N<n\le 2N}f_H^2(n)\ll N, 
$$
\par
\noindent
giving in all a new formula involving our correlations (from previous sections).  Thus, 
$$%\label{symDis} (36)
J_{f_H\, ,\,\sgn}(N,H)=2\sum_{h>0}W_H(h)\left(C_{f_H}(2N,h)-C_{f_H}(N,h)\right) + O(NH+H^3). 
\leqno{(36)}
$$
\par
\noindent
We leave, as an exercise, to prove from the definition 
$$
J_{f_H\, ,\,\sgn}(N,H)\gg NH^2, 
$$
\par
\noindent
say, the $f_H$ symmetry integral is \lq \lq trivial\rq \rq. However, assuming $f_H$ is in the first class, 
$$
C_{f_H}(N,h)=\sum_{\ell \ll N}\widehat{C_{f_H}}(N,\ell)c_{\ell}(h), 
\quad \hbox{\rm with}\quad 
\widehat{C_{f_H}}(N,\ell)\ll_{\varepsilon} {{N^{1+\varepsilon}}\over {\ell^2}}, 
$$
$$
C_{f_H}(2N,h)=\sum_{\ell \ll N}\widehat{C_{f_H}}(2N,\ell)c_{\ell}(h), 
\quad \hbox{\rm with}\quad 
\widehat{C_{f_H}}(2N,\ell)\ll_{\varepsilon} {{N^{1+\varepsilon}}\over {\ell^2}}, 
$$
\par				% PAGE 23 
\noindent
entailing 
$$
2\sum_{h>0}W_H(h)\left(C_{f_H}(2N,h)-C_{f_H}(N,h)\right)
\ll_{\varepsilon} \sum_{\ell \ll N}{{N^{1+\varepsilon}}\over {\ell^2}}\left|\sum_{h>0}W_H(h)c_{\ell}(h)\right|, 
$$
\par
\noindent
in which, defining as usual $e_q(n)\defineq e^{2\pi in/q}$ (the usual {\it additive characters} modulo $q$) 
$$
\sum_{h>0}W_H(h)c_{\ell}(h)={1\over 2}\sum_{h}W_H(h)c_{\ell}(h)-H\varphi(\ell) 
={1\over 2}\sum_{j\in \Z_{\ell}^*}\sum_{h}W_H(h)e_{\ell}(jh)-H\varphi(\ell) 
$$
\par
\noindent
using the $W_H(h)$ positive exponential sums [CL2] (with the additive characters orthogonality), 
$$
\sum_{h>0}W_H(h)c_{\ell}(h)\ll \sum_{j\le \ell}\sum_{h}W_H(h)e_{\ell}(jh)+H\varphi(\ell)
\ll \ell \sum_{h\equiv 0\bmod \ell}W_H(h)+H\varphi(\ell)
\ll \ell H  
$$
\par
\noindent
(uniformly in $\ell>1$, otherwise $c_1(h)=1$ gives $\sum_h W_H(h)=0$, compare [CL2]), so from $(36)$ we get 
$$
J_{f_H\, ,\,\sgn}(N,H)\ll_{\varepsilon} N^{1+\varepsilon}H+H^3, 
$$
\par
\noindent
say, in case $H\ll \sqrt{N}$, 
$$
NH^2\ll J_{f_H\, ,\,\sgn}(N,H)\ll_{\varepsilon} N^{1+\varepsilon}H, 
$$
\par
\noindent
a contradiction whenever  $N^{2\varepsilon}\ll H\ll \sqrt{N}$. 
\smallskip 
\par
Thus $f_H$ we built above is not in first class. 
Its shift-Ramanujan expansion is not regular and we guess that the real problem, 
here, is the dependence of our arithmetic function on the \lq \lq external\rq \rq  parameter $H$.
\bigskip
\bigskip
\par
\noindent {\bf Acknowledgements.}  We thank Biswajyoti Saha for his careful reading and corrections of an earlier version of this paper.

\bigskip
\bigskip
\bigskip

\par
\centerline{\stampatello References}
\bigskip

\item{\bf [Car]}
R. Carmichael, {Expansions of arithmetical functions in infinite series}, {\sl Proc. London Math.
Society}, {\bf 34} (1932), 1--26. 

\item{\bf [CojM]}
A.C. Cojocaru and M.R. Murty, {\sl An introduction to sieve methods and
their applications},  London Mathematical Society Student Texts, {\bf 66},
Cambridge University Press, Cambridge, 2006.

\item{\bf [CL1]} 
G. Coppola and M. Laporta, {Generations of correlation averages}, 
{\sl J.of Numbers}, {\bf 2014} (2014), Article ID 140840, 13 pages. 

\item{\bf [CL2]} 
G. Coppola and M. Laporta, {Symmetry and short interval mean-squares}, 
{\tt arxiv.org:}1312.5701, submitted 

\item{\bf [CL3]} 
G. Coppola and M. Laporta, {Sieve functions in arithmetic bands}, {\sl Hardy-Ramanujan J.}, {\bf 39} (2016), 21--37. 

\item{\bf [CMS]}
G. Coppola, M.R. Murty and B. Saha,
{On the error term in a Parseval type formula in the theory of
Ramanujan expansions II}, {\sl J. Number Theory}, {\bf 160} (2016), 700--715.

\item{\bf [CMS2]}
G. Coppola, M.R. Murty and B. Saha,
{Finite Ramanujan expansions and shifted convolution sums of arithmetical functions}, 
{\sl J. Number Theory}, {\bf 174} (2017), 78--92. 

\item{\bf [CSal]}
G. Coppola and S. Salerno, {On the symmetry of the divisor function
in almost all short intervals}, {\sl Acta Arith.}, {\bf 113} (2004), no. 2, 189--201. 

\item{\bf [GMP]}  
H.G. Gadiyar, M.R. Murty and R. Padma,
{Ramanujan - Fourier series and a theorem of Ingham},
{\sl Indian J. Pure Appl. Math.}, {\bf 45} (2014), no. 5, 691--706.

\item{\bf [HLi]}
G.H. Hardy, and J.E. Littlewood, {Some problems of \lq \lq Partitio numerorum\rq \rq. III: On the expression of a number as a sum of primes}
{\sl Acta Math.}, {\bf 44} (1923), 1--70. 
				% PAGE 24 
\item{\bf [MoV]}
H. Montgomery, R. Vaughan, {\sl Multiplicative Number Theory. Classical Theory}, 
Cambridge studies in advanced mathematics, {\bf 97}, Cambridge University Press, Cambridge, 2007. 

\item{\bf [M]}
M.R. Murty, {Ramanujan series for arithmetical functions},
{\sl Hardy-Ramanujan J.}, {\bf 36} (2013), 21--33. 

\item{\bf [MS]}
M.R. Murty and B. Saha,
{On the error term in a Parseval type formula in the theory of
Ramanujan expansions}, {\sl J. Number Theory}, {\bf 156} (2015), 125--134.

\item{\bf [R]}
S. Ramanujan, {On certain trigonometrical sums and their
applications in the theory of numbers}, {\sl Trans. Cambridge Philos.
Soc.}, {\bf 22} (1918), no. 13, 259--276.

\item{\bf [SchSpi]}
W. Schwarz and J. Spilker, 
{\sl Arithmetical functions,  (An introduction to elementary and analytic properties of arithmetic functions and to some of their almost-periodic properties).} London Mathematical Society Lecture Note Series, {\bf 184}, Cambridge University Press, Cambridge, 1994. 

\item{\bf [W]}
A. Wintner, {\sl Eratosthenian averages}, Waverly Press, Baltimore, MD, 1943. 

\bigskip
\bigskip
\bigskip

\par
\leftline{\tt Giovanni Coppola \hfill M. Ram Murty}
\leftline{\tt Universit\`a degli Studi di Napoli, \hfill Department of Mathematics,}
\leftline{\tt Complesso di Monte S. Angelo-Via Cinthia \hfill Queen's University,}
\leftline{\tt 80126 Napoli (NA), Italy \hfill Kingston, Ontario,}
\leftline{\tt www.giovannicoppola.name \hfill K7L 3N6, Canada}
\leftline{\tt giovanni.coppola@unina.it \hfill murty@mast.queensu.ca}

\bye